\documentstyle[10pt]{article}
\begin{document}

\baselineskip 15pt
\parindent=1em
\hsize=12.3 cm \textwidth=12.3 cm
\vsize=18.5 cm \textwidth=18.5 cm

\def\supt{{\rm supt}}
\def\dom{{\rm dom}}
\def\bfone{{\bf 1}}
\def\Gen{{\rm Gen}}
\def\rk{{\rm rk}}
\def\cali{{\cal I}}
\def\calj{{\cal J}}
\def\Spec{{\rm Spec}}

\title{Understanding preservation theorems: Chapter VI of {\em Proper and Improper Forcing}}
\author{
Chaz Schlindwein \\
Department of Mathematics and Computing \\
Lander University \\
Greenwood, South Carolina 29649 USA\\
{\tt cschlind@lander.edu}}

\maketitle

\def\forces{\mathbin{\parallel\mkern-9mu-}}
\def\notforces{\,\nobreak\not\nobreak\!\nobreak\forces}

\def\restr{\,\hbox{\vrule height8pt width.4pt depth0pt
   \vrule height7.75pt width0.3pt depth-7.5pt\hskip-.2pt
   \vrule height7.5pt width0.3pt depth-7.25pt\hskip-.2pt
   \vrule height7.25pt width0.3pt depth-7pt\hskip-.2pt
   \vrule height7pt width0.3pt depth-6.75pt\hskip-.2pt
   \vrule height6.75pt width0.3pt depth-6.5pt\hskip-.2pt
   \vrule height6.5pt width0.3pt depth-6.25pt\hskip-.2pt
   \vrule height6.25pt width0.3pt depth-6pt\hskip-.2pt
   \vrule height6pt width0.3pt depth-5.75pt\hskip-.2pt
   \vrule height5.75pt width0.3pt depth-5.5pt\hskip-.2pt
   \vrule height5.5pt width0.3pt depth-5.25pt}\,}

   \centerline{{\bf Abstract}}

   We present an exposition of Section VI.1 and most of Section VI.2 from Shelah's book {\em Proper and Improper Forcing}.  These sections offer proofs of the preservation under countable support iteration of proper forcing of various properties, including proofs that $\omega^\omega$-bounding, the Sacks property, the Laver property, and the $P$-point property are preserved by countable support iteration of proper forcing.

\eject

\section{Introduction}

This paper is an exposition of some preservation theorems, due to Shelah [12, Chapter VI], for countable support iterations of proper forcing. 
These include the preservation of the ${}^\omega\omega$-bounding property, the Sacks and Laver properties, the $P$-point property, and some others.  Generalizations to revised countable support iterations of semi-proper forcings or even certain non-semi-proper forcings are given in [13, Chapter VI] but we do not address these more general iterations.  The results of [12, Section VI.2] overlap the results of [2] and [3], but the methods are dissimilar.

This is the third in a sequence of expository papers covering parts of Shelah's book, {\it Proper and Improper Forcing.}  The earlier papers were [11], which covers sections 2 through 8 of [12, Chapter XI] and [9], which covers sections 2 and 3 of [12, Chapter XV].  Other papers by the author generalize certain other results in [12]; in no instance were we content to quote a result of Shelah without supplying a proof.  Thus, [5] may be read, in part, as an exposition of [12, Sections V.6, IX.2, and IX.4]; [6] is, in part, an exposition of [12, Section V.8 and Theorem III.8.5];
and [7] includes as a special case an alternative proof of [12, Theorem III.8.6].  Also, [5] answers
[12, Question IX.4.9(1)];  [6] answers a question implicit in [12, Section IX.4]; [10] answers another such question and also may be read, in part, as an exposition of the results of Eisworth and Shelah [1] that weaken the assumption
$``\alpha$-proper for every $\alpha<\omega_1$'' used in [12, Section V.6]. Lastly, [8, Section 2] corrects some minor errors in [5].

\section{Preservation of properness}

The fact that properness is preserved under countable support iterations was proved by Shelah in 1978.  The proof of this fact 
 is the basis of all preservation theorems for countable support iterations.

\proclaim Theorem 2.1 (Proper Iteration Lemma, Shelah). Suppose $\langle P_\eta\,
\colon\allowbreak\eta\leq\kappa\rangle$
is a countable support forcing iteration based on\/
$\langle\dot Q_\eta\,\colon\allowbreak\eta<\kappa\rangle$ and
for every $\eta<\kappa$ we have that\/ ${\bf 1}\forces_{P_\eta}``\dot Q_\eta$ {\rm is
proper.''} Suppose also that\/ $\alpha<\kappa$ and\/
$\lambda$ is a sufficiently large regular cardinal and\/ $N$ is a
countable elementary submodel of\/ $H_\lambda$ and\/
$\{P_\kappa,\alpha\}\in N$ and\/
$p\in P_\alpha$ is\/ $N$-generic and\/ {\rm
$p\forces``q\in  P_{\alpha,\kappa}\cap
N[G_{P_\alpha}]$.''}
Then there is $r\in P_\kappa$ such that
$r$ is $N$-generic and $r\restr\alpha=p$ and\/ {\rm
$p\forces``r\restr[\alpha,\kappa)\leq
q$.''}

Proof.  The proof proceeds by induction, so suppose that
the Theorem holds for all iterations of length less than $\kappa$.  Fix $\lambda$ a sufficiently large regular cardinal, and fix $N$ a countable elementary substructure of $H_\lambda$ such that $P_\kappa\in N$ and fix also $\alpha\in\kappa\cap N$ and $p\in P_\alpha$ and
a $P_\alpha$-name $q$ such that $p$ is $N$-generic and
$p\forces``q\in P_{\alpha,\kappa}\cap N[G_{P_\alpha}]$.''  

Case 1.  $\kappa=\beta+1$ for some $\beta$.

Because $\beta\in N$  we may use the induction hypothesis to fix $p'\in P_\beta$ such that
$p'\restr\alpha=p$ and $p'$ is $N$-generic and
$p\forces``p'\leq q\restr\beta$.''  We have that
$p'\forces``q(\beta)\in N[G_{P_\beta}]$.'' Take $r\in P_\kappa$ such that
$r\restr\beta=p'$ and 

\medskip

{\centerline{$p'\forces``r(\beta)\leq q(\beta)$ and $r(\beta)$ is
$N[G_{P_\beta}]$-generic for $Q_\beta$.''}}

\medskip

  Then $r$ is $N$-generic and we are done with the successor case.

Case 2.  $\kappa$ is a limit ordinal.

Let $\beta={\rm sup}(\kappa\cap N)$, and fix $\langle\alpha_n\,\colon\allowbreak
n\in\omega\rangle$ an increasing sequence from $\kappa\cap N$ cofinal in $\beta$ such that
$\alpha_0=\alpha$.  Let $\langle\sigma_n\,\colon\allowbreak n\in\omega\rangle$ enumerate all the $P_\kappa$ names $\sigma\in N$ such that ${\bf 1}\forces``\sigma$ is an ordinal.''

Using the induction hypothesis, build a sequence $\langle \langle p_n,q_n\rangle\,\colon\allowbreak n\in\omega\rangle$ such that $p_0=p$ and $q_0=q$ and for each $n\in\omega$ we have all of the following:

(1) $p_n\in P_{\alpha_n}$ and $p_n$ is $N$-generic  and
$p_{n+1}\restr\alpha_n=p_n$.

(2)  $p_n\forces``q_n\in P_{\alpha_n,\kappa}\cap N[G_{P_{\alpha_n}}]$ and
if $n>0$ then $q_n\leq q_{n-1}\restr[\alpha_n,\kappa)$ and $q_n\forces`\sigma_{n-1}\in
N[G_{P_{\alpha_n}}]$.'\thinspace''

(3) $p_n\forces``p_{n+1}\restr[\alpha_n,\alpha_{n+1})\leq q_n\restr\alpha_{n+1}$.''

Define $r\in P_\kappa$ such that $(\forall n\in\omega)\allowbreak
(r\restr\alpha_n=p_n)$ and ${\rm supp}(r)\subseteq\beta$.  
To see that $r$ is $N$-generic, suppose that $\sigma\in N$ is a $P_\kappa$-name for an ordinal. Fix $n$ such that
$\sigma = \sigma_n$. 
Because $p_{n+1}$ is $N$-generic, we have 

\medskip

\centerline{$p_{n+1}\forces``{\rm supp}(q_{n+1})\subseteq\kappa\cap N[G_{P_{\alpha_{n+1}}}]=
\kappa\cap N$,''}

\medskip

\noindent whence it is clear that 

\medskip

\centerline{$p_{n+1}\forces``r\restr[\alpha_{n+1},\kappa)\leq 
q_{n+1}$.'' }

\medskip

We have 

\medskip

\centerline{$p_{n+1}\forces``q_{n+1}\forces`\sigma\in Ord\cap N[G_{P_{\alpha_{n+1}}}]=Ord 
\cap N
$,'\thinspace''}

\medskip

\noindent where $Ord$ is the class of all ordinals. 
Thus $r\forces``\sigma\in N$.''
We conclude that $r$ is $N$-generic, and the Theorem is established.

\proclaim Corollary 2.2 (Fundamental Theorem of Proper Forcing, Shelah).
Suppose $\langle P_\eta\,
\colon\allowbreak\eta\leq\kappa\rangle$
is a countable support forcing iteration based on\/
$\langle Q_\eta\,\colon\allowbreak\eta<\kappa\rangle$ and
for every $\eta<\kappa$ we have that\/ ${\bf 1}\forces_{P_\eta}`` Q_\eta$ {\rm is
proper.''}  Then $P_\kappa$ is proper.

Proof: Take $\alpha = 0$ in the Proper Iteration Lemma.

\section{Preservation of proper plus $\omega^\omega$-bounding}

In this section we recount Shelah's proof of the preservation of ``proper plus $\omega^\omega$-bounding.''  This is a special case of [12, Theorem VI.1.12].
Another treatment of this result can be found in [2] and [3], using different methods.

\proclaim Lemma 3.1.  Suppose $\langle P_\eta\,\colon\allowbreak
\eta\leq\kappa\rangle$ is a countable support iteration based on $\langle Q_\eta\,\colon\allowbreak \eta<\kappa\rangle$ and each $Q_\eta$ is proper in
$V[G_{P_\eta}]$.
Suppose ${\rm cf}(\kappa)=\omega$ and
 $\langle\alpha_n\,\colon\allowbreak n\in\omega\rangle$ is an increasing sequence of ordinals cofinal in $\kappa$ with $\alpha_0=0$.
Suppose also that $f$ is a $P_\kappa$-name for an element
of ${}^\omega\omega$, and suppose $p\in P_\kappa$.
Then there are
$\langle p_n\,\colon\allowbreak n\in\omega\rangle$ and $\langle f_n\,\colon\allowbreak n\in\omega\rangle$  such that
$p_0\leq p$ and for every $n\in\omega$ we have that each of the following holds:

(1) For all $k\leq n$ we have ${\bf 1}\forces_{P_{\alpha_n}}``p_0\restr[\alpha_n,\kappa)
\forces`f(k)=f_n(k)$,'\thinspace'' and

(2) $f_n$ is a $P_{\alpha_n}$-name for an element of ${}^\omega\omega$, and

(3) $p_0\restr\alpha_n\forces``p_0\restr[\alpha_n,\alpha_{n+1})\forces`f_n(k)=f_{n+1}(k)$ for every $k\leq n+1$,'\thinspace'' and

(4) $p_{n+1}\leq p_n$, and

(5) whenever $k\leq m<\omega$ we have $p_0\restr\alpha_n\forces``p_m\restr[
\alpha_n,\alpha_{n+1})\forces`f_n(k)=f_{n+1}(k)$.'\thinspace''

\medskip

Proof:  
Fix $\lambda$ a sufficiently large regular cardinal and let $N$ be a countable
elementary substructure of $H_\lambda$ containing $P_\kappa$ and
$\langle\alpha_n\,\colon\allowbreak n\in\omega\rangle$ and $f$ and $p$.
Build $\langle p'_n\,\colon\allowbreak n\in\omega\rangle$ and
 $\langle\sigma_n\,\colon\allowbreak n\in\omega\rangle$ such that
$p'_0=p$ and each of the following holds:

(1) $p'_{n+1}\restr\alpha_n=p'_n\restr\alpha_n$.

(2) ${\bf 1}\forces_{P_{\alpha_n}}``\sigma_n\in\omega$ and
$p'_{n+1}\restr[\alpha_n,\kappa)
\forces`f(n)=\sigma_n$.'\thinspace''

(3) $p'_n\restr\alpha_n\forces``p'_{n+1}\restr[\alpha_n,\kappa)\leq 
p'_n\restr[\alpha_n,\kappa)$.''

(4) $p'_{n}\restr\alpha_n$ is $N$-generic.

(5) $p'_{n}\restr\alpha_n\forces``p'_{n}\restr[\alpha_n,\kappa)\in 
P_{\alpha_n,\kappa}\cap N[G_{P_{\alpha_n}}]$.''

(6)  $p'_{n}\in P_\kappa$.

Notice that (6) does not follow from the fact that $p'_{n}\restr\alpha_n\in
P_{\alpha_n}$ and $p'_{n}\restr\alpha_n\forces``p'_{n}\restr[\alpha,\kappa)\in
P_{\alpha,\kappa}$,'' but it does follow from (4) and (5).

Let $q_0=\bigcup\{p'_n\restr\alpha_n\,\colon\allowbreak n\in\omega\}$.  

At this point we define 
$f_n(k)=\sigma_k$ for $k\leq n$.  We have yet to define $f_n(k)$ for $k>n$.  Notice 
that we cannot set $f_n(k)=\sigma_k$ for $k>n$ because in $V[G_{P_{\alpha_n}}]$ 
we have that $\sigma_k$ is not an integer, but only a name.

Claim. For all $k\leq n$ we have $q_0\restr\alpha_n\forces``q_0\restr[\alpha_n,\kappa)
\forces`f(k)=f_n(k)$.'\thinspace''

Proof. Obvious.

Fix $\lambda'$ a suffciently large regular cardinal and $M$ a 
countable elementary substructure of
$H_{\lambda'}$ containing $N$ and $q_0$.

Build $\langle p^n_0\,\colon\allowbreak n\in\omega\rangle$
and $\langle\tau_n\,\colon\allowbreak n\in\omega\rangle$
 such that $p^0_0=q_0$ and each of the following holds:

(1) $p_0^{n+1}\restr\alpha_n=p^n_0\restr\alpha_n$.

(2) $p_0^{n+1}\leq p_0^n$.

(3)  ${\bf 1}\forces_{P_{\alpha_n}}``\tau_n\in\omega$ and $
p_0^{n+1}\restr[\alpha_n,\alpha_{n+1})\forces`\tau_n=f_{n+1}(n+1)$.'\thinspace''

(4) $p^n_0\restr\alpha_n$ is $M$-generic.

(5) $p^n_0\restr\alpha_n\forces``p^{n+1}_0\restr[\alpha_n,\kappa)\in P_{\alpha_n,\kappa}\cap
M[G_{P_{\alpha_n}}]$.''

Notice that (1), (4), and (5) imply that $p_0^{n+1}\in P_\kappa$; this is the reason the
structure $M$ is needed.

There is no difficulty in doing this.  At this point, we define
$f_n(n+1)=\tau_n$ for every $n$ and we let $p_0=\bigcup\{p^n_0\restr\alpha_n\,\colon\allowbreak n\in\omega\}$.  

At this point the following parts of the Lemma are exemplified:

(1) For all $k\leq n$ we have $p_0\restr\alpha_n\forces``p_0\restr[\alpha_n,\kappa)
\forces`f(k)=f_n(k)$.'\thinspace''

(2) $f_n\restr(n+2)$ is a $P_{\alpha_n}$-name for an element of ${}^{n+2}\omega$.

(3) ${\bf 1}\forces_{P_{\alpha_n}}``p_0\restr[\alpha_n,\alpha_{n+1})\forces`f_n(k)=f_{n+1}(k)$ for every $k\leq n+1$.'\thinspace''

Choose $\lambda^*$ a sufficiently large regular cardinal.
We build $\langle p_n\,\colon\allowbreak n\in\omega\rangle$ 
and $\langle M_n\,\colon\allowbreak n\in\omega\rangle$ by recursion on $n\in\omega$.
Let $M_0$ be a countable elementary substructure of $H_{\lambda^*}$ containing
$M$ and $p_0$. 

Fix $n$, and suppose $p_n$ and $M_n$ have been defined.

For each $i<n$
let $q^i_n$ and $\xi^i_{n+1}$ be chosen such that

\medskip

\centerline{${\bf 1}\forces_{P_{\alpha_i}}``\xi^i_{n+1}\in\omega$ and
$q^i_n\in P_{\alpha_i,\alpha_{i+1}}\cap M_n[G_{P_{\alpha_i}}]$ and}

\centerline{$q^i_n\leq p_n\restr[\alpha_i,\alpha_{i+1})$ and
$q^i_{n}\forces`f_{i+1}(n+1)=\xi^i_{n+1}$.'\thinspace''}

\medskip

Build $\langle r^i_n\,\colon\allowbreak i\leq n\rangle$ such that for each $i\leq n$ we have the following.

(1) $r^i_n\in P_{\alpha_i}$.

(2) $r^i_n\leq p_n\restr\alpha_i$.

(3) $r^i_n$ is $M_n$-generic.

(4) If $i<n$ then $r^{i+1}_n\restr\alpha_i=r^i_n$.

(5) If $i<n$ then $r^i_n\forces``r^{i+1}_n\restr[\alpha_i,\alpha_{i+1})\leq q^i_n$.''

Then take $p_{n+1}$ such that 

\medskip

\centerline{$p_{n+1}\restr\alpha_{n}=\bigcup\{r^i_n\,\colon\allowbreak i\leq n\}$ and
$p_{n+1}\restr[\alpha_{n},\kappa)=p_{n}\restr[\alpha_{n},\kappa)$.}

\medskip

 Let $M_{n+1}$ be a countable 
elementary substructure of $H_{\lambda^*}$ containing $M_n$ and $p_{n+1}$.

This completes the recursive construction.

We set $f_i(k)=\xi^i_k$ whenever $i+1<k$.

This completes the proof of the Lemma.

\proclaim Definition 3.2.  For $f$ and $g$ in ${}^\omega\omega$ we say
$f\leq g$ iff $(\forall n\in\omega)\allowbreak(f(n)\leq g(n))$.
We say that $P$ is ${}^\omega\omega$-bounding iff 
$V[G_P]\models``(\forall f\in{}^\omega\omega)\allowbreak(\exists g\in {}^\omega\omega\cap V)
\allowbreak(f\leq g)$.''

\proclaim Theorem 3.3.  Suppose $\langle P_\eta\,\colon\allowbreak\eta\leq\kappa\rangle$ is a countable support iteration based on $\langle Q_\eta\,\colon\allowbreak\eta<\kappa\rangle$ and suppose\/ {\rm $(\forall\eta<\kappa)\allowbreak({\bf 1}\forces_{P_\eta}``Q_\eta$ 
is proper and ${}^\omega\omega$-bounding'').} Suppose $f$ is a $P_\kappa$-name for an element 
of\/ ${}^\omega\omega$.
Then whenever $\lambda$ is a sufficiently large regular cardinal and\/ $N$ is a countable 
elementary substructure of $H_\lambda$
and $\alpha<\kappa$ and $\{P_\kappa,\alpha,f\}\in N$ and $p\in P_\alpha$ and $p$ is $N$-generic  
then 
{\rm $p\forces``(\forall q\in P_{\alpha,\kappa}\cap N[G_{P_\alpha}])\allowbreak(\exists q^\#\leq q)\allowbreak(\exists h^\#\in{}^\omega\omega)\allowbreak
(q^\#\forces`f\leq h^\#$').''}

Proof: The proof proceeds by induction on $\kappa$.  We assume that $\lambda$, 
$N$, $\alpha$, $p$, and $f$ are as in the hypothesis of the Theorem.
Fix  $q$  a $P_\alpha$-name in $N$ such that ${\bf 1}\forces``q\in P_{\alpha,\kappa}$.''

Case 1. $\kappa=\beta+1$. 

 Because ${\bf 1}\forces_{P_\beta}``Q_\beta$ is ${}^\omega\omega$-bounding,''
we may take $q^*$ and $h^*$ to be $P_\beta$-names such that

\medskip

\centerline{
${\bf 1}\forces_{P_\beta}``q^*\leq q(\beta)$ and $h^*\in{}^\omega\omega$ and
$q^*\forces`f\leq h^*$.'\thinspace''}

\medskip

We may assume that the names $q^*$ and $h^*$ are in $N$.
By the induction hypothesis we may take $P_\alpha$-names
$\tilde q$  and $h$  such  that

\medskip

\centerline{
$p\forces``\tilde q\leq q\restr\beta$ and $h\in{}^\omega\omega$ and
$\tilde q\forces`h^*\leq h$.'\thinspace''}

\medskip

  Define $q'$ such that
$p\forces``q'=(\tilde q, q^*)\in P_{\alpha,\kappa}$.''
Clearly 

\medskip

\centerline{$p\forces``q'\leq q$ and $q'\forces`f\leq h$.'\thinspace''}

\medskip

This completes Case 1.

\medskip

Case 2.  ${\rm cf}(\kappa)>\omega$.

Because no $\omega$-sequences of ordinals can be added at limit stages of uncountable cofinality,
 we may take $\beta$ and $f'$ and $q'$ to be $P_\alpha$-names in $N$ such that

\medskip

\centerline{
${\bf 1}\forces``\alpha\leq\beta<\kappa$ and ${\bf 1}\forces_{P_{\alpha,\beta}}`f'\in{}^\omega\omega$' and
 $q'\leq q$ and}

\centerline{
$q'\forces_{P_{\alpha,\kappa}}`f'=f$.'\thinspace''}

\medskip

For every $\beta_0\in\kappa\cap N$ such that $\alpha\leq\beta_0$ 
let $q^*(\beta_0)$ and $h(\beta_0)$ be  $P_\alpha$-names in $N$ such that

\medskip

\centerline{
${\bf 1}\forces``$if $\beta=\beta_0$ and there is some $q^*\leq q'\restr\beta$ and some 
$h\in{}^\omega\omega$ such that }

\centerline{
$q^*\forces`f'\leq h$,' then $q^*(\beta_0)$ and $h(\beta_0)$ are witnesses thereto.''}

\medskip

Let $q^*$ and $h$ and $s$ be $P_\alpha$-names such that for every $\beta_0
\in \kappa\cap N$, if $\alpha\leq\beta_0$ then

\medskip

\centerline{
${\bf 1}\forces``$if $\beta=\beta_0$ then $q^*=q^*(\beta_0)$ and $h=h(\beta_0)$ and $s\in P_{\alpha,\kappa}$}

\centerline{ and $s\restr\beta=q^*$ and
$s\restr[\beta,\kappa)=q'\restr[\beta,\kappa)$.''}

\medskip

Claim 1: $p\forces``s\leq q$ and $h\in{}^\omega\omega$  and $s
\forces`f\leq h$.'\thinspace''

Proof:  Suppose $p'\leq p$. Fix $p^\#\leq p'$ and $\beta_0<\kappa$ such that
$p^\#\forces``\beta_0=\beta$.''  Because the name $\beta$ is in $N$ and $p^\#$ is $N$-generic,
we have that $\beta_0\in N$.  Notice by the induction hypothesis that
we have

\medskip

\centerline{$p\forces``$there is some $q^\#\leq q'\restr\beta_0$ and some 
$h^\#\in{}^\omega\omega$}

{\centerline{such that $q^\#\forces`
f'\leq h^\#$.'\thinspace''}

\medskip

Hence 

\medskip

\centerline{$p^\#\forces``q^*=q^*(\beta_0)\leq q'\restr\beta$ and $h=h(\beta_0)\in{}^\omega\omega$ and}

\centerline{
$q^*\forces`f'\leq h$ and
$q'\restr[\beta,\kappa)\forces``f'=f$.''\thinspace'\thinspace''}

\medskip

\noindent  Therefore
$p^\#\forces``s\forces`f\leq h$.'\thinspace''

We conclude that
$p\forces``s\forces`f\leq h$.'\thinspace''
Claim 1 is established.

This completes Case 2.

\medskip

Case 3. ${\rm cf}(\kappa)=\omega$.

Let $\langle\alpha_n\,\colon\allowbreak n\in\omega\rangle$ be an increasing 
sequence from $\kappa\cap N$ cofinal in $\kappa$ such that 
$\alpha_0=\alpha$.

Let $\langle g_j\,\colon\allowbreak j<\omega\rangle$ list every $P_\alpha$-name
$g\in N$ such that
${\bf 1}\forces_{P_\alpha}`` g\in{}^\omega\omega$.''

Fix $\langle (p_n,f_n)\,\colon\allowbreak n\in\omega\rangle$ as in
Lemma 3.1 (applied in $V[G_{P_\alpha}])$.  That is, ${\bf 1}\forces``p_0\leq q$'' and 
 for every $n\in\omega$ we have that each of the following holds in $V[G_{P_\alpha}]$:

(0) $p_n\in P_{\alpha,\kappa}$.

(1) For every $k\leq n$ we have $p_0\restr\alpha_n\forces``p_0\restr[\alpha_n,\kappa)
\forces`f(k)=f_n(k)$.'\thinspace''

(2) $f_n$ is a $P_{\alpha_n}$-name for an element of ${}^\omega\omega$.

(3) $p_0\restr\alpha_n\forces``p_0\restr[\alpha_n,\alpha_{n+1})\forces`f_n(k)=f_{n+1}(k)$ for every $k\leq n+1$.'\thinspace''

(4) $p_{n+1}\leq p_n$.

(5) Whenever $k\leq m<\omega$ we have $p_0\restr\alpha_n\forces``p_m\restr[
\alpha_n,\alpha_{n+1})\forces`f_n(k)=f_{n+1}(k)$.'\thinspace''

We may assume that for every $n\in\omega$ the $P_\alpha$-names $p_n$ and $f_n$
are in $N$, and, furthermore, the sequence $\langle\langle 
p_n,f_n\rangle\,\colon\allowbreak n\in\omega\rangle$ is in $ N$.

In $V[G_{P_\alpha}]$, define $\langle g^n\,\colon\allowbreak n\in \omega\rangle$ by
\medskip

{\centerline{$g^n(k)={\rm max}\{f_0(k),\allowbreak{\rm max}\{g_j(k)\,\colon\allowbreak
j\leq n\}\}$.}}

\medskip

 Also in $V[G_{P_\alpha}]$ define $g\in{}^\omega\omega$ such that

\medskip

{\centerline{$g(k)=g^k(k)$ for all $k\in\omega$.}}

\medskip

Claim 2.  Suppose $\alpha\leq\beta\leq\gamma<\kappa$  and
suppose $f'$ is a $P_\gamma$-name for an element of ${}^\omega\omega$. Then 

\medskip

\centerline{${\bf 1}\forces_{P_\beta}``V[G_{P_\alpha}]\models`(\forall q\in P_{\beta,\gamma})(\exists q'\leq q)(\exists h'\in {}^\omega\omega)
(q'\forces``f'\leq h'$'').'\thinspace''}

\medskip

Proof: Given $r_1\in P_\alpha$ and a $P_\alpha$-name $r_2$ for an element of $P_{\alpha,\beta}$ and a $P_\beta$-name
$q$ for an element of $P_{\beta,\gamma}\cap V[G_{P_\alpha}]$, choose $\lambda'$ a sufficiently large regular cardinal and
$N'$ a countable elementary substructure of $H_{\lambda'}$ containing $\{r_1,r_2,q,P_\kappa, \alpha,\beta,\gamma,f'\}$.
Choose $r_1'\leq r_1$ such that $r_1'$ is $N'$-generic.
By the overall induction hypothesis (i.e., because $\gamma<\kappa$) we have

\medskip

\centerline{$r_1'\forces``(\exists s\leq (r_2,q))(\exists h'\in{}^\omega\omega)(s\forces`f'\leq h'$').''}

\medskip

Consequently we may fix $s$ and $h'$ such that

\medskip

\centerline{$(r_1',s\restr\beta)\forces``V[G_{P_\alpha}]\models`s\restr[\beta,\gamma)\leq q$ and
 $s\restr[\beta,\gamma)\forces``f'\leq h'$.''\thinspace'\thinspace''}

\medskip

The Claim is established.

Claim 3.  We may be build $\langle r_n\,\colon\allowbreak n\in\omega\rangle$ such that
$r_0=p$  and for every $n\in\omega$ we have that the following hold:

(1) $r_n\in P_{\alpha_n}$ is $N$-generic.

(2) $r_{n+1}\restr\alpha_n=r_n$.

(3) $r_n\forces``f_n\leq g$.''

(4) $p\forces``r_n\restr[\alpha,\alpha_n)\leq p_0\restr[\alpha,\alpha_n)$.''

Proof: Work by induction on $n$.  For $n=0$ there is nothing to prove, so assume that
$n>0$ and suppose we have $r_n$.

Fix $P_{\alpha_n}$-names $F_0$ and $F_2$
such that ${\bf 1}\forces``$if there are functions $F_0'$ and $F_2'$
such that $F'_0\in V[G_{P_\alpha}]$ maps $P_{\alpha_n,\alpha_{n+1}}$ 
into ${}^\omega\omega$ and
$F'_2\in V[G_{P_\alpha}]$ maps $P_{\alpha_n,\alpha_{n+1}}$ into $P_{\alpha_n,\alpha_{n+1}}$ and
for every $q'\in P_{\alpha_n,\alpha_{n+1}}\cap V[G_{P_\alpha}]$ we have
$F'_2(q')\leq q'$ and $F'_2(q')\forces`f_{n+1}\leq F'_0(q')$', then $F_0$ and $F_2$ are witnesses
to this.''

We may assume that the names $F_0$ and $F_2$ are in $N$.

By Claim 2 we have 

\medskip

(*) $r_n\forces``F_0\in V[G_{P_\alpha}]$ maps $P_{\alpha_n,\alpha_{n+1}}$ 
into ${}^\omega\omega$ and
$F_2\in V[G_{P_\alpha}]$ maps $P_{\alpha_n,\alpha_{n+1}}$ into $P_{\alpha_n,\alpha_{n+1}}$ and
for every $q'\in P_{\alpha_n,\alpha_{n+1}}\cap V[G_{P_\alpha}]$ we have
$F_2(q')\leq q'$ and $F_2(q')\forces`f_{n+1}\leq F_0(q')$.'\thinspace''

\medskip

In $V[G_{P_{\alpha_n}}]$, define $ g_n^*$ by
$(\forall i\in\omega)(g_n^*(i)= 
{\rm max}\{F_0(p_m\restr[\alpha_n,\alpha_{n+1}))(i)\,\colon\allowbreak m\leq i\})$.

We may assume the name $g_n^*$ is in $N$.

Notice that we have

\medskip

\centerline{$r_n\forces``g_n^*\in N[G_{P_{\alpha_n}}]\cap V[G_{P_\alpha}] = N[G_{P_\alpha}]
$.''}

\medskip

Therefore we may choose a $P_{\alpha_n}$-name $k$ such that $r_n\forces``g_n^*=g_k$'' (in our notation, 
we suppress the fact that $k$ depends on $n$).

Subclaim 1: $r_n\forces``F_2(p_k\restr[\alpha_n,\alpha_{n+1}))
\forces`f_{n+1}\leq g$.'\thinspace''

Proof: 
For $i\geq k$ we have 

\medskip

{\centerline{$r_n\forces``F_2(p_k\restr[\alpha_n,\alpha_{n+1}))
\forces`f_{n+1}(i)\leq F_0(p_k\restr[\alpha_n,\alpha_{n+1}))(i)$}}

{\centerline{$
\leq g_n^*(i)=g_k(i)\leq g^i(i)=g(i)$.'\thinspace''}}

\medskip

\noindent The first inequality is by (*), 
the second inequality is by the definition of $g_n^*$ 
along with the fact that $i\geq k$, the equality is by the 
definition of $k$, the next inequality is by the definition 
of $g^i$ along with the fact that $i\geq k$, and the last 
equality is by the definition of $g$.

For $i<k$, we have 

\medskip

\centerline{$r_n\forces``p_k\restr[\alpha_n,\alpha_{n+1})
\forces`f_{n+1}(i)=f_n(i)\leq g(i)$.'\thinspace''}

\medskip

\noindent The equality is by the choice of $\langle (f_m, p_m)\,\colon m\in\omega\rangle$ (see Lemma 3.1),
and the inequality is by the induction hypothesis that Claim 3 holds for integers less than or equal to $n$.

Because $r_n\forces``F_2(p_k\restr[\alpha_n,\alpha_{n+1}))\leq p_k\restr
[\alpha_n,\alpha_{n+1})$,'' we have that the Subclaim is established.

Using the Proper Iteration Lemma, choose $r_{n+1}\in P_{\alpha_{n+1}}$ such that
$r_{n+1}$ is $N$-generic and $r_{n+1}\restr\alpha_n=r_n$ and

\medskip

\centerline{
$r_n\forces``r_{n+1}\restr[\alpha_n,\alpha_{n+1})\leq F_2(p_k\restr[\alpha_n,\alpha_{n+1}))$.''}

\medskip

This completes the proof of Claim 3.

Let $r'=\bigcup\{r_n\,\colon\allowbreak n\in\omega\}$.  We have that

\medskip

\centerline{
$p\forces``r'\restr[\alpha,\kappa)\leq q$ and 
$r'\restr[\alpha,\kappa)\forces`f\leq g$.'\thinspace''}

\medskip

The Theorem is established.

\proclaim Corollary 3.4.  Suppose $\langle P_\eta\,\colon\allowbreak\eta\leq\kappa\rangle$ is a countable support iteration based on $\langle Q_\eta\,\colon\allowbreak\eta<\kappa\rangle$ and suppose\/ {\rm $(\forall\eta<\kappa)\allowbreak({\bf 1}\forces_{P_\eta}``Q_\eta$ 
is proper and ${}^\omega\omega$-bounding'').}  Then $P_\kappa$ is $\omega^\omega$-bounding.

Proof.  Take $\alpha = 0$ in Theorem 3.3.

\section{The Sacks property}

In this section we present Shelah's proof of the preservation of ``proper plus Sacks property'' under countable support iteration.  The proof is a special case of  [12, Theorem VI.1.12].  

\proclaim Definition 4.1.  For $x$ and $y$ in ${}^\omega(\omega-\{0\})$,
 we say that $x\ll y$ iff $(\forall n\in\omega)\allowbreak(x(n)\leq y(n))$ and
\[ \lim_{n\rightarrow\infty}y(n)/x(n)=\infty \]
In particular for $x\in{}^\omega(\omega-\{0\})$ we have $1\ll x$ iff $(\forall k\in\omega)\allowbreak
(\exists n\in\omega)\allowbreak(\forall m\geq n)\allowbreak(x(m)>k)$.

\proclaim Definition 4.2.  
For $T\subseteq{}^{<\omega}\omega$ a tree and $x\in{}^\omega(\omega-\{0\})$, 
we say that $T$ is an $x$-sized tree iff for every $n\in\omega$ we have that the 
cardinality of $T\cap{}^n\omega$ is at most $x(n)$.

\proclaim Definition 4.3.  For $T\subseteq{}^{<\omega}\omega$ we set $[T]$
equal to the set of all $f\in{}^\omega\omega$ such that
every initial segment of $f$ is in $T$.  That is, $[T]$ is the set of infinite branches of $T$.

\proclaim Definition 4.4.  
A poset $P$ has the Sacks property iff whenever $x\in{}^\omega(\omega-\{0\})$ and 
$1\ll x$
 then we have

{\centerline{${\bf 1}\forces_P``(\forall f\in{}^\omega\omega)\allowbreak(\exists H\in V)\allowbreak(H$ is an $x$-sized tree and $f\in[ H])$.''}}

\proclaim Definition 4.5.  Suppose $n\in\omega$. We say that $t$ is an $n$-tree iff $t\subseteq{}^{\leq n}\omega$ and $t$ is closed under initial segments and $t$ is non-empty and for every $\eta\in t$ there is $\nu\in t$ such that $\nu $ extends $\eta$ and ${\rm lh}(\nu)=n$.

\proclaim Lemma 4.6.  Suppose $P$ has the Sacks property and $x$ and  $z$ are elements of ${}^\omega(\omega-\{0\})$ and $x\ll z$. Then we have

{\centerline{${\bf 1}\forces_P``(\forall T)($if $T$ is an $x$-sized tree then}}
{\centerline{$(\exists H\in V)\allowbreak(H$ is a $z$-sized tree and $T\subseteq H))$.''}}

\medskip

Proof: Work in $V[G_P]$. For every $n\in\omega$ let 

\medskip

\centerline{${\cal T}_n(x)=\{t\subseteq{}^{\leq n}\omega\,\colon\allowbreak
t$ is an $n$-tree}

\centerline{and $(\forall i\leq n)\allowbreak(\vert t\cap {}^i\omega\vert\leq x(i))\}$.}

\medskip

Let

\medskip

\centerline{${\cal T}(x)=\bigcup\{{\cal T}_n(x)\,\colon\allowbreak n\in\omega\}$.}

\medskip

 Under the natural order, 
${\cal T}(x)$ is isomorphic to ${}^{<\omega}\omega$.

Define $\zeta\in[{\cal T}(x)]$ by setting $\zeta(n)=T\cap {}^{\leq n}\omega$ for all $n\in\omega$.

Define $y\in {}^\omega(\omega-\{0\})$ by setting $y(n)$ equal to the greatest integer less than or equal to $z(n)/x(n)$ for every $n\in\omega$.
Clearly $1\ll y$, so we may choose a $y$-sized tree $H'\subseteq {\cal T}(x)$ such that $\zeta\in[H']$ and
$H'\in V$.

Let $H^*=\bigcup H'$ and let 

\medskip

\centerline{$H=\{\eta\in H^*\,\colon\allowbreak(\forall n\in\omega)\allowbreak (\exists\nu\in{}^n\omega\cap H^*)\allowbreak
(\nu$ is comparable with $\eta)\}$.}

\medskip

We have that $H$ is a $z$-sized tree and $H\in V$ and $T\subseteq H$.

The Lemma is established.

\proclaim Lemma 4.7. Suppose $n^*$ is an integer and suppose $y$ and $z$ are elements of ${}^\omega\omega$
and $y\ll z$. Suppose $P$ is a forcing such that\/ {\rm
$V[G_P]\models``$for every countable $X\subseteq V$ there is a countable $Y\in V$ such that $X\subseteq Y$.''}
Suppose in $V[G_P]$ we have a sequence $\langle T_n\,\colon\allowbreak n\in\omega\rangle$
such that for every $n$ we have $T_n\in V $ is a $y$-sized tree.
Then in $V[G_P]$ there is a $z$-sized tree $T^*\in V$ and an increasing sequence of integers $\langle k(n)\,\colon n\in\omega\rangle$
such that $k(0)=0 $ and $k(1)\geq n^*$ and $(\forall n>0)\allowbreak(k(n)>n)$ and for every $\eta\in {}^{<\omega}\omega
$ we have

\medskip

\centerline{$(\forall n\in\omega)\allowbreak(\exists i<n)(\eta\restr k({n+1})\in T_{k(i)})$ implies $\eta\in T^*$.}

\medskip

Proof: Fix $x\in {}^\omega(\omega-\{0\})$ such that $y\ll x\ll z$. Fix $\langle x_n\,\colon\allowbreak
n\in\omega\rangle$ a sequence of elements of ${}^\omega(\omega-\{0\})$ such that
$(\forall n\in\omega)(y\ll x_n\ll x_{n+1}\ll x$.

Work in $V[G_P]$. Let $b\in V$ be a countable set of $y$-sized trees such that $\{T_n\,\colon\allowbreak n\in\omega\}\subseteq b$.
Let $\langle S_n\,\colon\allowbreak n\in\omega\rangle\in V$ enumerate $b$ with infinitely many repetitions with $S_0=T_0$.

Define $h\in{}^\omega\omega$ by setting $h(0)=0$ and for every $n>0$ setting
$h(n)$ equal to the least $m>h(n-1)$ such that $T_n=S_m$.

For each $n\in\omega$ set $k(n)$ equal to the least $k\geq n^*$ such that

\medskip

\centerline{$(\forall j\geq k)(2x_n(j)\leq x_{n+1}(j)$ and $(n+2)y(j)\leq z(j))$.}

\medskip

Build $\langle S'_n\,\colon n\in\omega\rangle$ by setting $S'_0=S_0$ and for every $n\in\omega$ let

\medskip

\centerline{$S'_{n+1}=\{\rho\in S_n\,\colon\allowbreak\rho\restr k(n)\in S'_n\}\cup S'_n$.}

\medskip

Claim 1. For all $n\in\omega$ we have that $S'_n$ is an $x_n$-sized tree.

Proof: By induction on $n$. Clearly $S'_0$ is an $x_0$-sized tree. For every $t<k(n)$ we have that
$\vert S'_{n+1}\cap{}^t\omega\vert\leq \vert S'_n\cap{}^t\omega\vert\leq x_n(t)\leq x_{n+1}(t)$.
For every $t\geq k(n)$ we have $\vert S'_{n+1}\cap{}^t\omega\vert\leq
\vert S'_n\cap{}^t\omega\vert+\vert S_n\cap{}^t\omega\vert\leq x_n(t)+y(t)\leq x_{n+1}(t)$. 
The Claim is established.

Let $T^*=\{\eta\in{}^{<\omega}\omega\,\colon\allowbreak (\forall n>0)(\exists i<n)\allowbreak(\eta\restr k(n)\in S'_{k(i)})\}$.

Claim 2. $T^*$ is a $z$-saized tree.

Proof.  Given $t\geq k(1)$,  choose $n\in\omega$ such that $k(n)\leq t<k({n+1})$.
We have 

\medskip

\centerline{$ T^*\cap {}^t\omega=\{\eta\in{}^t\omega\,\colon\allowbreak(\forall j\leq n+1)\allowbreak
(\exists i\leq j)\allowbreak(\eta\restr k(j)\in T_{k(i)})\}$}

\medskip

\noindent and so 

\medskip

\centerline{$\vert T^*\cap {}^t\omega\vert\leq\Sigma_{i\leq n+1}\vert T_{k(i)}\cap {}^t\omega
\vert\leq(n+2)y(t)
\leq z(t)$.}

\medskip

For $t<k(1)$ we have $T^*\cap{}^t\omega=T_0\cap{}^t\omega$, so $\vert T^*\cap{}^t\omega\vert\leq y(t)\leq z(t)$.

The Claim is established

Build $\langle n'_i\,\colon\allowbreak i\in\omega\rangle$ an increasing sequence of integers such that
$n'_0=0$ and $n'_1>k(1)$ and for every $i\in\omega$ we have

(A) $ h(n'_{i})<n'_{i+1}$ and

(B) $k(n'_i)<n'_{i+1}$ and

(C) $(\exists t)(n'_i<k(t)<n'_{i+1})$.

For every $i\in\omega$ let $m_i=h(n'_{4i+4})$.

Fix $\eta\in{}^{<\omega}\omega$ such that

(D) $(\forall i>0)(\exists j<0)(\eta\restr m_{i+1}\in T_{m_j})$.

To establish the Lemma, it suffices to show 

(E) $(\forall i>0)(\exists j<i)(\eta\restr k(i)\in S'_{k(j)})$,

\noindent since this implies $\eta\in T^*$.

Claim 3. $(\forall i>0)(\exists j<i)(\eta\restr n'_{i+1}\in S'_{n'_j})$.

We prove this by induction on $i$.

Case 1. $i<9$.

We have $n_{i+1}'\leq n'_8\leq h(n'_8)=m_1$ and we have $\eta\restr m_1\in T_0$.
Therefore $\eta\restr n'_{n+1}\in T_0=S_0=S'_0$.

Case 2. $i\geq 9$,

Fix $i^*$ such that $4i^*+2\leq i<4i^*+6$.

By (D) we may fix $j^*<i^*$ such that $\eta\restr m_{i^*+1}\in T^*_{m_{j^*}}$.

Because $n'_{i+1}\leq h(n'_{4i^*+8})\leq m_{i^*+1}$, we have $\eta\restr n'_{i+1}\in T_{m_{j^*}}$.

If $m_{j^*}=0$, we are done, so assume otherwise.

Subclaim 1. Suppose $\rho\in T_{m_{j^*}}$ and $\rho\restr k(n'_{4j^*+4})\in S'_{n'_{4j^*+4}}$. Then
$\rho\in S'_{n'_{4j^*+5}}$.

Proof: We have $T_{m_{j^*}}=T_{h(n'_{4j^*+4})}=S_{n'_{4j^*+4}}$. Therefore
$\rho\in S'_{n'_{4j^*+4}+1}\subseteq S'_{n'_{4j^*+5}}$.

The Subclaim is established.

Subclaim 2. $k(n'_{4j^*+4})\leq m_{i^*+1}$.

Proof: $m_{i^*+1}=h(n'_{4i^*+8})>n'_{4i^*+7}\geq n'_{4j^*+11}\geq n'_{4j^*+6}>k(n'_{4j^*+4})$. The Subclaim is established.

Let $\rho=\eta\restr m_{i^*+1}$.

By the choice of $j^*$ we have 

(F) $\rho\in T_{m_{j^*}}$.

By Subclaim 2 we have 

(G) $\rho\restr k(n'_{4j^*+4})=\eta\restr k(n'_{4j^*+4})$.

Subclaim 3. $\eta\restr n'_{4j^*+5}\in S'_{n'_{4j^*+4}}$.

Proof: Because $4j^*+4\leq i$ we may use the induction hypothesis of the proof of Claim 3. The Subclaim is established.

By (B), (G), and Subclaim 3,  we have $\rho\restr k(n'_{4j^*+4})\in S'_{n'_{4j^*+4}}$.

Therefore by Subclaim 1, we have $\rho\in S'_{n'_{4j^*+5}}$.

To complete the proof of the Lemma, suppose $i>0$. We must show that there is $t<i$ such that
$\eta\restr k(i)\in S'_{k(t)}$.

Case 1: $k({i-1})<n'_0$.

By (C) we have $n'_1\geq k({i})$. By Claim 3 we have $\eta\restr n'_2\in S_0$.  Hence $\eta\restr n'_1\in S_0$. Hence
$\eta\restr k(i)\in S_0$.

Case 2:  $n'_0\leq k({i-1})$.

By (A) we know that there is at most one element of $\{n'_j\,\colon\allowbreak j\in\omega\}$ strictly between
$k({i-1})$ and $k(i)$. Hence we may
fix $j>0$ such that $n'_{j-1}\leq k({i-1})<k({i})\leq n'_{j+1}$.  If $\eta\restr n'_{j+1}\in S_0$ then
$\eta\restr k({i})\in S_0$ and we are done, so assume otherwise.  By Claim 3 we may fix $m<j$ such that $\eta\restr n'_{j+1}\in S'_{n'_m}$.
We have $\eta\restr k(i)\in S'_{n'_m}\subseteq S'_{n'_{j-1}}\subseteq S'_{k_{i-1}}$ and again we are done.

The Lemma is established.

\proclaim Lemma 4.8.  Suppose  $x\in{}^\omega(\omega-\{0\})$ 
and $z\in{}^\omega(\omega-\{0\})$ 
and  $x\ll z$.  
Suppose that for every $n\in\omega$ we have
that
$T_n$ is an $x$-sized tree. 
Suppose $T$ is an $x$-sized tree.
Then there is a $z$-sized tree $T^*\supseteq T$ and a sequence of integers
$\langle m_i\,\colon\allowbreak i\in\omega\rangle$
such that
for every $\eta\in T$ and $i\in\omega$ and every $\nu\in T_{m_i}$ extending $\eta$, if 
${\rm length}(\eta)\geq m_i$  then  $\nu\in T^*$.

Proof. Choose $y\in{}^\omega(\omega-\{0\})$ such that $x\ll y\ll z$.
Fix $n^*\in\omega$ such that $(\forall n\geq n^*)\allowbreak(2x(n)\leq y(n))$.
 For every $n\geq n^*$ define $T'_n=\{\eta\in{}^{<\omega}\omega\,\colon\allowbreak
\eta\in T$ or $\eta\in T_n$ and $\eta\restr n\in T\}$.
For every $n<n^*$ let $T'_n=T$.

For every $n\in\omega$ we have that $T'_n$ is a $y$-sized tree.

By Lemma 4.7 we may choose $T^*$ a $z$-sized tree and
$\langle k_i\,\colon\allowbreak i\in\omega\rangle$ an increasing sequence of integers
such that $k_0=0$ and $k_1\geq n^*$ and

\medskip

{\centerline{$(\forall \eta\in{}^{<\omega}\omega)\allowbreak((\forall n\in\omega)\allowbreak
(\exists i\leq n)\allowbreak(\eta\restr k_n\in T'_{k_i})$ implies $\eta\in T^*)$.}}

\medskip

Clearly $T\subseteq T^*$.

For every $i\in\omega$ set $m_i=k_{i+1}$.

Now suppose that $\eta\in T$ and $i\in\omega$ and
${\rm length}(\eta)\geq m_i$ and
 $\nu$ extends $\eta$
and $\nu\in T_{m_i}$.  We show $\nu\in T^*$.

Because $\nu$ extends an element of $T$ of length at least $m_i$, we have that $\nu\in T'_{m_i}$.
Choose $h\in[T'_{m_i}]$ extending $\nu$.  It suffices to show that $h\in[T^*]$.
Therefore it suffices to show that for every $n\in\omega$ we have

\begin{itemize}

\item[$(*)_n$]  $(\exists j\leq n)(h\restr k_n\in T'_{k_j})$.

\end{itemize}

Fix $n\in\omega$.

Case 1: $i<n$.

Because $h\in[T'_{k_{i+1}}]$ we have that $h\in[T'_{k_n}]$. Therefore
$h\restr k_n\in T'_{k_n}$ and we have that $(*)_n$ holds.  

Case 2: $n\leq i$.

We have $h\restr k_n=\eta\restr k_n\in T=T'_0$.  Therefore $(*)_n$ holds.

The Lemma is established.

\proclaim Lemma 4.9.  Suppose  
$y\in{}^\omega(\omega-\{0\})$ and $z\in{}^\omega(\omega-\{0\})$, and suppose
$\langle x_n\,\colon\allowbreak n\in\omega\rangle$ is a sequence of elements
of ${}^\omega(\omega-\{0\})$ such that $(\forall n\in\omega)\allowbreak
(x_n\ll x_{n+1}\ll y\ll z)$.  Suppose $T$ is an $x_0$-sized tree.
Suppose for every $n\in\omega$, we have $x_n^*\in{}^\omega(\omega-\{0\})$ and 
$x_n^*\ll x_n\ll x^*_{n+1}$, and we have
$\langle x_{n,j}\,\colon\allowbreak n\in\omega$, $j\in\omega\rangle$ is a
sequence of elements of ${}^\omega(\omega-\{0\})$ such that
for every  $j\in\omega$ we have
$x_n\ll x_{n,j}\ll x_{n,j+1}\ll x_{n+1}^*$.
Suppose  $\langle T_{n,j}\,\colon\allowbreak n\in\omega$, $j\in\omega\rangle$ is
a sequence such that for every $n\in\omega$ and
 $j\in\omega$ we have that $T_{n,j}$ is an
$x_{n,j}$-sized tree.  Then there are $\langle T^n\,\colon\allowbreak n\in\omega\rangle$ and
$T^*$ such that $T^*$ is a $z$-sized tree and $T\subseteq T^*$ and for every
$n\in\omega$ we have

(i) $T^n\subseteq T^{n+1}$ and $T^n$ is an $x_n$-sized tree, and

(ii) for every forcing notion $P$ we have that in $V[G_P]$  for every $j\in\omega$ and every $g\in[T_{n,j}]\cap V[G_P]$ 
there is $k\in\omega$ such that for every $\eta\in T_{n,j}$ extending $g\restr k$, if 
 $\eta\restr k\in T^n\cap T^*$ then
$\eta\in T^{n+1}\cap T^*$.

\medskip

Proof:    Let $T^0=T$. 
Given $T^n$, build $\langle T'_{n,j}\,\colon\allowbreak j\in\omega\rangle$
as follows.  Let $T'_{n,0}=T^n$.  Given
$T'_{n,j}$ take $m(n,j)\in\omega$ such that

\medskip

\centerline{
$(\forall t\geq m(n,j))\allowbreak(2x_{n,j}(t)\leq x_{n,j+1}(t))$.}

\medskip

Let $T'_{n,j+1}=\{\eta\in T_{n,j}\colon\allowbreak
\eta\restr m(n,j)\in T'_{n,j}\}\cup T'_{n,j}$.

Claim 1.  Whenever $i\leq j<\omega$ we have $T'_{n,i}\subseteq T'_{n,j}$.

Proof.  Clear.

Claim 2. Suppose $T^n$ is an $x_n$-sized tree.  Then $(\forall j\in\omega)\allowbreak(T'_{n,j}$
is an $x_{n,j}$-sized tree).

Proof:  It is clear that $T'_{n,0}$ is an $x_{n,0}$-sized tree.  
Assume that $T'_{n,j}$ is an $x_{n,j}$-sized tree.  Fix $t\in\omega$. 

Case 1: 
 $t<m(n,j)$.

We have that

\medskip

\centerline{ $T'_{n,j+1}\cap {}^t\omega=
T'_{n,j}\cap{}^t\omega$}

\medskip

\noindent and so 

\medskip

\centerline{$\vert T'_{n,j+1}\cap {}^t\omega\vert\leq x_{n,j}(t)\leq
x_{n,j+1}(t)$.}

\medskip

Case 2: $t\geq m(n,j)$.  

We have 

\medskip

\centerline{$T'_{n,j+1}\cap{}^t\omega\subseteq(T'_{n,j}\cap{}^t\omega)\cup (T_{n,j}\cap
{}^t\omega)$.}

\medskip

\noindent  Therefore
we have 

\medskip

\centerline{$\vert T'_{n,j+1}\cap{}^t\omega\vert\leq 2x_{n,j}(t)\leq x_{n,j+1}(t)$.}

\medskip

  The
Claim is established.

For each $n\in\omega$, using Claim 2 and Lemma 4.7 we my find an increasing sequence of integers $\langle k_{n,j}\,\colon\allowbreak j\in\omega\rangle$ and $T^{n+1}$ such that $k_{n,0}=0$ and
$(\forall j>0)\allowbreak(k_{n,j}>j)$ and
if $T^n$ is an $x_n$-sized tree, then $T^{n+1}$ is
an $x_{n+1}$-sized tree such that
for all $\eta\in{}^{<\omega}\omega$, we have 

\medskip

{\centerline{$(\forall j\in\omega)\allowbreak
(\exists i\leq j)\allowbreak
(\eta\restr k_{n,j}\in T'_{n,k_{n,i}})$ implies $\eta\in T^{n+1}.$}}

\medskip

This completes the construcion of $\langle T^n\,\colon\allowbreak n\in\omega\rangle$ and
$\langle T_{n,j}\,\colon\allowbreak j\in\omega$, $n\in\omega\rangle$.

Applying mathematical induction, we have that each $T^n$ is in fact an $x_n$-sized tree.

Claim 3. $T^n\subseteq T^{n+1}$ for every $n\in\omega$.

Proof: By Claim 1 we have that $T^n\subseteq T'_{n,i}$ for every $i\in\omega$.
By the definition of $T^{n+1}$ we have that

\medskip

\centerline{
$T^{n+1}\supseteq\bigcap\{T'_{n,k_{n,i}}\,\colon\allowbreak i\in\omega\}\supseteq
\bigcap\{T'_{n,i}\,\colon\allowbreak i\in\omega\}\supseteq T^n$.}

\medskip

The Claim is established.

Applying Lemma 4.7 again we obtain an increasing sequence of integers $\langle k_n\,\colon\allowbreak n\in\omega\rangle$ and a
$z$-sized tree $T^*$ such that $(\forall n\in\omega)\allowbreak
(n<k_n)$ and for every $\eta\in{}^{<\omega}\omega$,
we have that

\medskip

{\centerline{$(\forall n\in\omega)(\exists i\leq n)(
\eta\restr k_n\in T^{k_i})$ implies $\eta\in T^*$.}}

\medskip

Notice that $T^0\subseteq\bigcap\{T^n\,\colon\allowbreak n\in\omega\}\subseteq T^*$.

Now we verify that $\langle T^n\,\colon\allowbreak n\in\omega\rangle$ and $T^*$ satisfy the
remaining conclusions of the Lemma.  Accordingly, fix a forcing notion $P$ and work in $V[G_P]$. Fix
$n\in\omega$ and $j\in\omega$ and
$g\in[T_{n,j}]$.
Let 

\medskip

\centerline{$k={\rm max}(k_n,\allowbreak{\rm max}\{k_{n,j'}\,\colon\allowbreak j'\leq j\},
\allowbreak{\rm max}\{m(n,j')\,\colon\allowbreak
j'\leq j\})$.}

\medskip

  Fix $\eta\in T_{n,j}$ extending $g\restr k$ and assume that
 $\eta\restr k\in T^n\cap T^*$.

Claim 4. $\eta\in T^{n+1}$.

Proof:  It suffices to show 

\medskip

\centerline{$(\forall j'\in\omega)\allowbreak
(\exists i\leq j')\allowbreak(\eta\restr k_{n,j'}\in T'_{n,k_{n,i}})$.}

\medskip

  Fix $j'\in\omega$
and let $i={\rm min}(j,j')$.

Case 1: $j'\leq j$. 

Because $k_{n,j'}\leq k$ we have that 
 $\eta\restr k_{n,j'}\in T^n\subseteq T'_{n,k_{n,i}}$, as required.

Case 2: $j<j'$.

It suffices to show that $\eta\restr k_{n,j'}\in T'_{n,k_{n,j}}$. 
Because $g\restr k=\eta\restr k\in T^n$ and $m(n,j)\leq k$, we have that
 $g\restr m(n,j)\in T^n\subseteq T'_{n,j}$.  Because
we have $\eta\in
T_{n,j}$ and $\eta\restr m(n,j)=g\restr m(n,j)\in T'_{n,j}$, we know
 by the definition of
$T'_{n,j+1}$ and Claim~1 that
$\eta\in T'_{n,j+1}\subseteq T'_{n,k_{n,j}}$.

Claim 4 is established.

Claim 5. $\eta\in T^*$.

Proof: It suffices to show $(\forall i\in\omega)\allowbreak(\exists i'\leq i)
\allowbreak(\eta\restr k_i\in T^{k_{i'}})$.  Towards this end, fix $i\in\omega$.

Case 1:  $i\leq n$.

Because $\eta\restr k\in T^*$ and $\eta$ extends $g\restr k$,
 we have $g\restr k_i\in T^*$ and hence
we may take $i'\leq i$ such that $g\restr k_i\in T^{k_{i'}}$. But we also have
 $\eta\restr k_i=g\restr k_i$, so we finish Case~1.

Case 2: $n<i$.

We let $i'=i$. 
By Claim 4 we have $\eta\restr k_i\in T^{n+1}$, and by Claim 3 we have that
$T^{n+1}\subseteq T^{k_n}\subseteq T^{k_i}$.

Claim 5 and the Lemma are established.

\proclaim Theorem 4.10.  Suppose $\langle P_\eta\,\colon\allowbreak\eta\leq\kappa\rangle$ is a countable support iteration based on $\langle Q_\eta\,\colon\allowbreak\eta<\kappa\rangle$ and suppose\/ {\rm $(\forall\eta<\kappa)\allowbreak({\bf 1}\forces_{P_\eta}``Q_\eta$ 
is proper and has the Sacks property'').} 
Suppose $\lambda$ is a sufficiently large regular cardinal 
and $\alpha<\kappa$ 
and  $x$ and $z$ are $P_\alpha$-names and $T$ is a $P_\kappa$-name 
and\/ {\rm ${\bf 1}\forces_{P_\alpha}``x\in
{}^\omega(\omega-\{0\})$ and $z\in{}^\omega(\omega-\{0\})$ and
$ x\ll z$''} and ${\bf 1}\forces_{P_\kappa}``T$ is an $x$-sized tree.''
Suppose $N$ is a countable elementary submodel of $H_\lambda$ and $\{P_\kappa, \alpha,  x, z, T
\}\in N$. 
Suppose
 $p\in P_\alpha$ and $p$ is $N$-generic.
Then {\rm
$p\forces``(\forall q\in P_{\alpha,\kappa}\cap N[P_{G_{P_\alpha}}])\allowbreak(\exists q'\leq q)\allowbreak(\exists H)\allowbreak(H$
is a $z$-sized tree and $q'\forces`T\subseteq H$').''}

Proof: The proof proceeds by induction on $\kappa$.  We assume that 
  $\lambda$, 
$N$, $\alpha$, $p$,  $x$, $z$, and $T$ are as in the hypothesis of the Theorem. Fix a $P_\alpha$-name $q$
in $N$ such that ${\bf 1}\forces``q\in P_{\alpha,\kappa}$.''

Case 1. $\kappa=\beta+1$.

Fix $y$ a $P_\alpha$-name in $N$ such that ${\bf 1}\forces``x\ll y\ll z$.''

Using Lemma 4.6, we may choose $\tilde q$ and $H'$ such that 

\medskip

\centerline{
${\bf 1}\forces_{P_\beta}``\tilde q\leq q(\beta)$ and $H'$ is a $y$-sized tree and
$\tilde q\forces`T\subseteq H'$.'\thinspace''}

\medskip

We may assume that the names $\tilde q$ and $H'$ are elements of $N$.
Using the induction hypothesis we get  $P_\alpha$-names $q^*$ and
 $H$  such that

\medskip

\centerline{
$p\forces``q^*\leq q\restr\beta$ and $H$ is an $x$-sized tree and
$q^*\forces`H'\subseteq H$.'\thinspace''}

\medskip

We have that $p\forces``(q^*,\tilde q)\leq q'$'' and Case 1 is established.

Case 2.  ${\rm cf}(\kappa)>\omega$.

Because no $\omega$-sequences of ordinals can be added at limit stages of uncountable cofinality, we may take $\beta$ and $T'$ and $q'$ to be $P_\alpha$-names in $N$ such that

\medskip

\centerline{
${\bf 1}\forces``\alpha\leq\beta<\kappa$ and $T'$ is a
$P_{\alpha,\beta}$-name and $q'\leq q$ and}

\centerline{
${\bf 1}\forces_{P_{\alpha,\beta}}`T'$ is an $x$-sized tree' and
$q'\forces_{P_{\alpha,\kappa}}`T'=T$.'\thinspace''}

\medskip

For every $\beta_0\in\kappa\cap N$ such that $\alpha\leq\beta_0$ 
let $\tilde q(\beta_0)$ and $H(\beta_0)$  be  $P_\alpha$-names in $N$ such that

\medskip

\centerline{
${\bf 1}\forces``$if $\beta=\beta_0$ and there is some $\tilde q\leq q'\restr\beta$ and some 
$H^*$}

\centerline{ such that $H^*$ is a $z$-sized tree and
$\tilde q\forces`T'\subseteq H^*$,'}

\centerline{
then $q^*(\beta_0)$ and $H(\beta_0)$  are witnesses thereto.''}

\medskip

Let $ q^*$ and $H$  and $s$ be $P_\alpha$-names such that for every $\beta_0
\in \kappa\cap N$, if $\alpha\leq\beta_0$, then

\medskip

\centerline{
${\bf 1}\forces``$if $\beta=\beta_0$ then $ q^*=q^*(\beta_0)$ and $H=H(\beta_0)$
 and 
$s\in P_{\alpha,\kappa}$ and}

\centerline{ $s\restr\beta=q^*$ and
$s\restr[\beta,\kappa)=q'\restr[\beta,\kappa)$.''}

\medskip

Claim 1: $p\forces``s\leq q$ and $s
\forces`T\subseteq H$.'\thinspace''

Proof:  Suppose $p'\leq p$. Fix $p^*\leq p'$ and $\beta_0<\kappa$ such that
$p^*\forces``\beta_0=\beta$.''  Because the name $\beta$ is in $N$ and $p^*$ is $N$-generic,
we have that $\beta_0\in N$.  Notice by the induction hypothesis
we have

\medskip

\centerline{$p\forces``$there is some $q^\#\leq q'\restr\beta_0$ and 
some $z$-sized tree
$H^\#$}

\centerline{such that $q^\#\forces`T'\subseteq H^\#$.'\thinspace''}

\medskip

Hence 

\medskip

\centerline{$p^*\forces``q^*=q^*(\beta_0)\leq q'\restr\beta$ and $H=H(\beta_0)$ and
$H$ is a $z$-sized tree and}

\centerline{$q^*\forces`T'\subseteq H$ and
$q'\restr[\beta,\kappa)\forces``T'=T$.''\thinspace'\thinspace''}

\medskip

  Therefore
$p^*\forces``s\forces`T\subseteq H$.'\thinspace''

Claim 1 is established.
This completes Case 2.

Case 3. ${\rm cf}(\kappa)=\omega$.

Let $\langle\alpha_n\,\colon\allowbreak n\in\omega\rangle$ be an increasing 
sequence from $\kappa\cap N$ cofinal in $\kappa$ such that 
$\alpha_0=\alpha$.

Define $y\in{}^\omega(\omega-\{0\})$ by letting $y(n)$ be the greatest integer less than or equal to
$z(n)/x(n)$ for every $n\in\omega$.  

In $V[G_{P_\alpha}]$, let $\langle {\cal T}_n(x)\,\colon\allowbreak n\in\omega\rangle$ and ${\cal T}(x)$ be as in the proof of Lemma 4.6. 
Also in $V[G_{P_\alpha}]$,
fix an isomorphism of 
${}^{<\omega}\omega$ onto ${\cal T}(x)$ and implicitly fix a $P_\alpha$-name for the isomorphism that is an element of $N$.

Let $\zeta$ be a $P_\kappa$-name in $N$ such that ${\bf 1}\forces_{P_\kappa}``\zeta\in[{\cal T}(x)]$ and for every $n\in\omega$ we have
 $\zeta(n)=T\cap {}^{\leq n}\omega$.''

Fix
 $\langle (p_n,\zeta_n)\,\colon\allowbreak n\in\omega\rangle\in N$ as in Lemma 3.1 applied in $V[G_{P_\alpha}]$
(that is, the sequence of names is an element of $N$ but not necessarily their values). That is, we have 
${\bf 1}\forces``p_0\leq q$'' and 
 for every $n\in\omega$ we have that each of the following holds:

(0) $p_n$ is a $P_\alpha$-name for an element of $P_{\alpha,\kappa}$.

(1) For every $k\leq n$ we have ${\bf 1}\forces_{P_{\alpha_n}}``p_0\restr[\alpha_n,\kappa)
\forces`\zeta\restr k=\zeta_n\restr k$.'\thinspace''

(2) $\zeta_n$ is a $P_{\alpha_n}$-name for an element of $[{\cal T}(x)]$.

(3) ${\bf 1}\forces_{P_{\alpha_n}}``p_0\restr[\alpha_n,\alpha_{n+1})\forces`\zeta_n\restr k
=\zeta_{n+1}\restr k$ for every $k\leq n+1$.'\thinspace''

(4) ${\bf 1}\forces_{P_\alpha}``p_{n+1}\leq p_n$.''

(5) Whenever $k\leq m<\omega$ we have ${\bf 1}\forces_{P_{\alpha_n}}``p_m\restr[
\alpha_n,\alpha_{n+1})\forces`\zeta_n\restr k=\zeta_{n+1}\restr k$.'\thinspace''

Claim 2.  Suppose $\alpha\leq\beta\leq\gamma<\kappa$ and ${\bf 1}\forces_{P_\alpha}
``x'\ll z'$''  and
suppose $T'$ is a $P_\gamma$-name for an $x'$-sized tree. Then

\medskip

\centerline{${\bf 1}\forces_{P_\beta}``V[G_{P_\alpha}]\models`(\forall q\in P_{\beta,\gamma})(\exists q'\leq q)(\exists H)$}

\centerline{$(H $ is a $z'$-sized tree and $q'\forces``T'\subseteq H$'').'\thinspace''}

\medskip

Proof: Given $r_1\in P_\alpha$ and a $P_\alpha$-name $r_2$ for an element of $P_{\alpha,\beta}$ and a 
$P_\beta$-name $q$ for an element of $P_{\beta,\gamma}\cap V[G_{P_\alpha}]$, choose $\lambda'$ a sufficiently large
regular cardinal and $N'$ a countable elementary substructure of $H_{\lambda'}$ containing $\{r_1,r_2,q,P_\kappa,\alpha,\beta,\gamma,x',z', T'\}$.
Choose $r_1'\leq r_1$ such that $r_1'$ is $N'$-generic.  By the overall induction hypothesis (i.e., because $\gamma<\kappa$)
we may choose $s$ such that

\medskip

\centerline{$r'_1\forces`` s\leq(r_2,q)$ and $(\exists H)(H$ is a $z'$-sized tree and $s\forces`T'\subseteq H$').''}

\medskip

Consequently we may choose $H$ such that

\medskip

\centerline{$(r'_1,s\restr\beta)\forces``V[G_{P_\alpha}]\models`
s\restr[\beta,\gamma)\leq q$ and $H$ is a $z'$-sized tree and}

\centerline{$s\restr[\beta,\gamma)\forces``T'\subseteq H$'').'\thinspace''}

\medskip

The Claim is established.

In $V[G_{P_\alpha}]$ fix $y'\ll y$ such that $1\ll y'$.  We may assume that the name $y'$ is in $N$.

Let $\Omega=\{x'\in N\,\colon\allowbreak x'$ is a $P_\alpha$-name and 
${\bf 1}\forces`` 1\ll x'\ll y'$''$\}$. Let $\langle y_n\,\colon\allowbreak n\in\omega\rangle$ enumerate $\Omega$.  Build $\langle x^*_n\,\colon\allowbreak n\in\omega\rangle$ as follows.  Let $x^*_0=y_0$, and for each $n\in\omega$ choose $x^*_{n+1}\in \Omega$ such that
${\bf 1}\forces``x^*_n\ll x^*_{n+1}$ and $y_{n+1}\ll x^*_{n+1}$.''
Also build $\langle x_n\,\colon\allowbreak n\in\omega\rangle$ a sequence of elements of $\Omega$ such that for each $n\in\omega$ we have
$x^*_n\ll x_n\ll x^*_{n+1}$.

For each $n\in\omega$ let $\langle T_{n,j}\,\colon\allowbreak j\in\omega\rangle$ list all $P_{\alpha}$-names 
$T'\in N$ such that in $V[G_{P_\alpha}]$ there is some $y'\ll x^*_{n+1}$ such that $T'\subseteq{\cal T}(x)$ is a
$y'$-sized tree, and build $\langle x_{n,j}\,\colon\allowbreak j\in\omega\rangle$ a sequence of elements of
$\Omega$ such that for every $j\in\omega$ we have in $V[G_P]$ that
$x_n\ll x_{n,j}\ll x_{n,j+1}\ll x^*_{n+1}$ and $T_{n,j}\subseteq{\cal T}(x)$ is an $x_{n,j}$-sized tree.

Using Lemma 4.9, choose $T^*\subseteq{\cal T}(x)$  a $y$-sized tree  and
 $\langle T^n\,\colon\allowbreak n\in\omega\rangle\in V[G_{P_\alpha}]$ a sequence of subsets of ${\cal T}(x) $ such that
$T^*\in V[G_{P_\alpha}]$ and
$T^0\subseteq T^*$ and $\zeta_0\in [T^0]$ and for every $n\in\omega$ we have that
$T^n$ is an $x_n$-sized tree and $T^n\subseteq T^{n+1}$ and, in $V[G_{P_\kappa}]$, we have that
 for every $j\in\omega$ and every 
$g\in[T_{n,j}]$  there is $k\in\omega$ such that for
every $\eta\in T_{n,j}$ extending $g\restr k$, if
 $\eta\restr k\in T^n\cap T^*$ then
$\eta\in T^{n+1}\cap T^*$. We may assume the $P_\alpha$-names $T^*$ and $\langle T^n\,\colon\allowbreak n\in\omega\rangle$ are in $N$.

Note that the reason we worked in $V[G_{P_\kappa}]$ rather then in $V[G_{P_\alpha}]$ in the previous paragraph is because 
we wish to allow $g$ to range over $[T_{n,j}]$ with the brackets interpreted in $ V[G_{P_\kappa}]$ (i.e., $g$ need not be in
$V[G_{P_\alpha}])$.

Using Claim 2, for every $n\in\omega$, let $F_{n,0}$ and $F_{n,2}$ and $y^*_n$ be $P_{\alpha_n}$-names  such that

\medskip

(A)  
 ${\bf 1}\forces``F_{n,0}$ and
$F_{n,2}$ and $y^*_n$ are functions, all three of which are in $V[G_{P_\alpha}]$,

\centerline{and each of whose domains
 is equal to $P_{\alpha_n,\alpha_{n+1}}$, such that}

\centerline{$(\forall q'\in P_{\alpha_n,\alpha_{n+1}}\cap V[G_{P_\alpha}]
 )\allowbreak (F_{n,0}(q')\subseteq{\cal T}(x)$ is a $y_n^*(q')$-sized tree}

\centerline{and $y^*_n(q')\ll x^*_{n+1}$ and
$F_{n,2}(q')\leq q'$}

\centerline{and $F_{n,2}(q')\forces`\zeta_{n+1}\in[ F_{n,0}(q')]$').''}

\medskip

 We may assume that the names $F_{n,0}$ and
$F_{n,2}$ and $y^*_n$ are in $N$.

For each $n\in\omega$ we may, in $V[G_{P_{\alpha_n}}]$, use Lemma 4.8 to choose $y_n\ll x_{n+1}$
and
$\tilde T_n\subseteq{\cal T}(x)$ a $y_n$-sized tree and
$\langle k^n_i\,\colon\allowbreak i\in\omega\rangle$ an increasing sequence of
integers 
such that 
 $T^n\subseteq\tilde T_n$ and for every $\eta\in T^n$ and every $i\in\omega$ and every
$\nu\in F_{n,0}(p_{k^n_i}\restr[\alpha_n,\alpha_{n+1}))$,
 if 
${\rm length}(\eta)\geq k^n_i$ and 
$\nu$ extends $\eta$, then $\nu\in\tilde T_n$. 

\medskip

We may assume the $P_{\alpha_n}$-names $\tilde T_n$ and $\langle k^n_i\,\colon\allowbreak i\in\omega\rangle$ are in $N$.

Claim 3.  We may be build $\langle r_n\,\colon\allowbreak n\in\omega\rangle$ such that
$r_0=p$  and for every $n\in\omega$ we have that the following hold:

(1) $r_n\in P_{\alpha_n}$ is $N$-generic, and

(2) $r_{n+1}\restr\alpha_n=r_n$, and

(3) $r_n\forces``\zeta_n\in [T^n]\cap [T^*]$,'' and

(4) $p\forces``r_n\restr[\alpha,\alpha_n)\leq p_0\restr\alpha_n$.''

\medskip

Proof: By induction on $n$.  For $n=0$ we have nothing to prove. Suppose we have $r_n$.

  By (A) and the definition of $\tilde T_n$  we have that

\medskip

(B)
$r_n\forces``T^n\subseteq\tilde T_n$ ''

and 

(C) $r_n\forces``$for every $\eta\in T^n$ and every $i\in\omega$

\centerline{and every
$\nu\in F_{n,0}(p_{k^n_i}\restr[\alpha_n,\alpha_{n+1}))$,}

\centerline{
 if 
${\rm length}(\eta)\geq k^n_i$ and 
$\nu$ extends $\eta$ then $\nu\in\tilde T_n$.''}

\medskip

By (C) and the fact that, by the induction hypothesis, we know $r_n\forces``\zeta_n\in[T^n]$,'' we have that

(D) $r_n\forces``(\forall j\in\omega)(F_{n,2}(p_{k_j^n}\restr[\alpha_n,\alpha_{n+1}))\forces`(\forall\nu\in
F_{n,0}(p_{k^n_j}\restr[\alpha_n,\alpha_{n+1})))\allowbreak($if $\nu$ extends $\zeta_n\restr k^n_j$ then $\nu\in\tilde T_n))$'\thinspace.''

We have

\medskip

{\centerline{$r_n\forces``\tilde T_n\in N[G_{P_{\alpha_n}}]$.''}}

\medskip

We also have

\medskip

{\centerline{$r_n\forces``\tilde T_n\in V[G_{P_{\alpha}}]$.''}}

\medskip

Therefore, because $r_n$ is $N$-generic, we have

\medskip

{\centerline{$r_n\forces``\tilde T_n\in N[G_{P_\alpha}]$.''}}

\medskip

Therefore there is a $P_{\alpha_n}$-name $m$ such that

\medskip

{\centerline{$r_n\forces``\tilde T_n=T_{n,m}$.''}}

\medskip

Using this fact along with the fact that $\langle T^n\,\colon\allowbreak n\in\omega\rangle$ and $T^*$ were chosen as in the conclusion of  Lemma 4.9 and also using the fact that $r_n\forces``\zeta_n\in[T^n]\subseteq[\tilde T_n]$,'' we may choose
 $k$ to be a $P_{\alpha_n}$-name for an integer such that
\medskip

(E) $r_n\forces``(\forall\eta\in \tilde T_n)\allowbreak($if $\eta$ extends $\zeta_n\restr k$ and $\eta\restr k\in T^n\cap T^*$ 
then $\eta\in T^{n+1}\cap T^*)$.''

\medskip

Choose $j$  to be a $P_{\alpha_n}$-name for an integer such that
$r_n\forces``k^n_j\geq k$.''

Subclaim 1. 
$r_n\forces``F_{n,2}(p_{k^n_j}\restr[\alpha_n,\alpha_{n+1}))\forces
`\zeta_{n+1}\in [\tilde T_n]$.'\thinspace''

Proof.  It suffices to show

\medskip 

$r_n\forces``F_{n,2}(p_{k^n_j}\restr[\alpha_n,\alpha_{n+1}))\forces
`(\forall j'>j)\allowbreak(\zeta_{n+1}\restr k^n_{j'}\in\tilde T_n)$.'\thinspace''

\medskip

Fix $j'$ a $P_{\alpha_{n+1}}$-name for an integer such that 

\medskip

$r_n\forces``F_{n,2}(p_{k^n_j}\restr[\alpha_n,\alpha_{n+1}))\forces j'>j$.'\thinspace''

\medskip

We know by the induction hypothesis that $r_n\forces``\zeta_n\in[T^n]$.'' Therefore

(F) $r_n\forces``\zeta_n\restr k^n_j\in T^n$.''

By the definition of $\langle p_i\,\colon\allowbreak i\in\omega\rangle$, we have

\medskip

(G) $r_n\forces``p_{k^n_j}\restr[\alpha_n,\alpha_{n+1})\forces`\zeta_n\restr k^n_j=
 \zeta_{n+1}\restr k^n_j$.'\thinspace''

\medskip

By (A) we have

\medskip

(H) $r_n\forces``F_{n,2}(p_{k^n_j}\restr[\alpha_n,\alpha_{n+1}))\forces`\zeta_{n+1}\in
[F_{n,0}(p_{k^n_j}\restr[\alpha_n,\alpha_{n+1}))]$.'\thinspace''

\medskip

Combining (G), (H), the definition of $\tilde T_n$, we have that

\medskip

{\centerline{$r_n\forces``F_{n,2}(p_{k^n_j}\restr[\alpha_n,\alpha_{n+1}))\forces`\zeta_{n+1}\restr k^n_{j'}\in\tilde T_n$.'\thinspace''}}

The Subclaim is established.

Subclaim 2. $r_n\forces``F_{n,2}(p_{k^n_j}\restr[\alpha_n,\alpha_{n+1}))\forces`\zeta_{n+1}\in [T^{n+1}]\cap [T^*]$.'\thinspace''

Proof: By (E) we have

\medskip

(I) $r_n\forces``(\forall\eta\in \tilde T_n)\allowbreak($if $\eta$ extends $\zeta_n\restr k^n_j$ and $\eta\restr k^n_j\in T^n\cap T^*$

\noindent then
$\eta\in T^{n+1}\cap T^*$).''

\medskip

Work in $V[G_{P_{\alpha_n}}]$ with $r_n\in G_{P_{\alpha_n}}$.
Fix $\eta\in\tilde T_n$ and
 suppose 

\medskip

\centerline{$F_{n,2}(p_{k^n_j}\restr[\alpha_n,\alpha_{n+1}))\forces``
\eta$ is an initial segment of $\zeta_{n+1}$ with
${\rm lh}(\eta)\geq k^n_j$.''}

\medskip

To establish the Subclaim, it suffices to show

\medskip

(J) $F_{n,2}(p_{k^n_j}\restr\alpha_n,\alpha_{n+1}))\forces``\eta\in T^{n+1}\cap T^*$.''

\medskip

By the definition of $\langle p_i\,\colon\allowbreak i\in\omega\rangle$ we have

\medskip

{\centerline{$p_{k^n_j}\restr[\alpha_n,\alpha_{n+1})\forces``\eta\restr k^n_j= \zeta_{n+1}\restr k^n_j=\zeta_n\restr k^n_j$.''}}

\medskip

Hence by the fact that Claim 3 holds for the integer $n$ we have

\medskip

(K) $p_{k^n_j}\restr[\alpha_n,\alpha_{n+1}))\forces``\eta\restr k^n_j\in T^n\cap T^*$.''

\medskip

By Subclaim 1, (I), (K), and the fact that
$F_2(p_{k^n_j}\restr[\alpha_n,\alpha_{n+1}))\leq
p_{k^n_j}\restr[\alpha_n,\alpha_{n+1})$ we obtain (J).

Subclaim 2 is established.

To complete the induction establishing Claim 3, we use the Proper Iteration Lemma to take
$r_{n+1}\in P_{\alpha_{n+1}}$ such that $r_{n+1}\restr\alpha_n=r_n$ and
$r_{n+1}$ is $N$-generic and $r_n\forces``r_{n+1}\restr[\alpha_n,\alpha_{n+1})\leq
F_{n,2}(p_{k^n_j}\restr[\alpha_n,\alpha_{n+1}))$.''

Claim 3 is established.

Let $q'$ be a $P_\alpha$-name such that

\medskip

\centerline{$p\forces``q'=\bigcup\{r_n\restr[\alpha,\alpha_n)\,\colon\allowbreak n\in\omega\}$.''}

\medskip

In $V[G_{P_\alpha}]$, let $H^*=\bigcup T^*$ and let

\medskip

\centerline{$H=\{\nu\in H^*\,\colon\allowbreak(\forall n\in\omega)\allowbreak(\exists\eta\in H^*)\allowbreak
(\nu$ is comparable with $\eta)\}$.}

\medskip

 As in the proof of Lemma 4.6, we have that $H$ is a $z$-sized tree. By Claim 3 we have
that

\medskip

\centerline{$q'\forces``$for every $n\in\omega$ we have $\zeta_n\in[T^*]$ and $\zeta_n\restr n=\zeta\restr n$,}

\centerline{and therefore $\zeta\in[T^*]$,
and therefore $T\subseteq H$.''}

\medskip

The Theorem is established.

\proclaim Corollary 4.11.  Suppose $\langle P_\eta\,\colon\allowbreak\eta\leq\kappa\rangle$ is a countable support iteration based on $\langle Q_\eta\,\colon\allowbreak\eta<\kappa\rangle$ and suppose\/ {\rm $(\forall\eta<\kappa)\allowbreak({\bf 1}\forces_{P_\eta}``Q_\eta$ 
is proper and has the Sacks proeprty.'')}  Then $P_\kappa$ has the Sacks property.

Proof.  Take $\alpha = 0$ in Theorem 4.10.

\section{The Laver Property}

In this section, we present Shelah's proof that the Laver property is preserved by countable support iteration of proper forcing.

\proclaim Definition 5.1. Suppose $f\in {}^\omega(\omega-\{0\})$ and ${\bf 1}\ll f$.
We say that $T$ is an $f$-tree iff $T$ is a tree and $(\forall\eta\in T)\allowbreak(\forall n\in{\rm dom}(\eta))\allowbreak(\eta(n)<f(n))$.

\proclaim Definition 5.2.
We say that $P$ is $f$-preserving iff whenever $z$ is
in ${}^\omega(\omega-\{0\})$ and $1\ll z$ then 

${\bf 1}\forces_P``(\forall g\in{}^\omega\omega)(g\leq f$ implies there exists $H\in V$ such that $H$ is
a $z$-sized $f$-tree  and $g\in[ H])$.''

\medskip

\proclaim Definition 5.3. We say that $P$ has the Laver property iff for every $f\in{}^\omega(\omega-\{0\})$ such that ${\bf 1}\ll f$ we have that $P$ is $f$-preserving.

\proclaim Theorem 5.4. $P$ has the Sacks property iff $P$ has the Laver property and $P$ is ${}^\omega\omega$-bounding.

Proof:  We first assume that $P$ has the Sacks property and we show that $P$ is ${}^\omega\omega$-bounding.  Given $p\in P$ and a name $f$ such that $p\forces``f\in{}^\omega\omega$,'' take $z\in{}^\omega(\omega-\{0\})$ such that ${\bf 1}\ll z$ and use the fact that $P$ has the Sacks property to obtain $q\leq p$ and a $z$-sized tree $H$ such that
$q\forces``f\in[H]$.'' For every $n\in\omega$ let 

\medskip

\centerline{$g(n)={\rm max}\{\eta(n)\,\colon\allowbreak\eta\in H$ and ${\rm lh}(\eta)>n\}$.}

\medskip

 Then we have $q\forces``f\leq g$.'' This establishes the fact that $P$ is ${}^\omega\omega$-bounding.

It is clear that if $P$ has the Sacks property, then it has the Laver property.

Finally we assume that $P$ has the Laver property and is ${}^\omega\omega$-bounding, and we show that $P$ has the Sacks property.
So suppose that $p\in P$ and ${\bf 1}\ll z$ and $p\forces``g\in{}^\omega\omega$.''
Using the fact that $P$ is ${}^\omega\omega$-bounding, take $p'\leq p$ and  $f\in{}^\omega\omega$ such that
$p'\forces``g\leq f$.'' Using the fact that $P$ has the Laver property, take $q\leq p'$ and $H$ a
$P$-name such that

\medskip

\centerline{$q\forces``H$ is a $z$-sized $f$-tree and $g\in[H]$ and $H\in V$.''}

\medskip

The Theorem is established.

\proclaim Theorem 5.5. Suppose $\langle P_\eta\,\colon\allowbreak\eta\leq\kappa\rangle$ is a countable support iteration based on $\langle Q_\eta\,\colon\allowbreak\eta<\kappa\rangle$ and suppose\/ {\rm $(\forall\eta<\kappa)\allowbreak({\bf 1}\forces_{P_\eta}``Q_\eta$ 
is proper and has the Laver property.'')}  Then $P_\kappa$ has the Laver property.

Proof: Fix $f\in{}^\omega\omega$ such that ${\bf 1}\ll f$.  Repeat the proofs of Lemma 4.5 through Corollary 4.11 with ``tree'' replaced by ``$f$-tree.'' The Theorem is established.

\section{$(f,g)$-bounding}

In this section we establish the preservation of $(f,g)$-bounding forcing.  For an exact formulation, see
 Corollary 6.8 below.  This proof is due to Shelah, of course; see [12, Conclusion VI.2.11F].

\proclaim Definition 6.1. We say that $T$ is an $(f,g)$-corseted tree iff

(0) $T\subseteq{}^{<\omega}\omega$ is a tree, and

(1)  $f$ and $g$ 
are functions with domain $\omega$, and

(2) $(\forall n\in\omega)\allowbreak(f(n)\in\{r\in{\bf R}\,\colon
\allowbreak 1<r\}\cup\{\omega\})$, and

(3) $(\forall n\in\omega)\allowbreak(g(n)\in\{r\in {\bf R}\,\colon\allowbreak
1<r\}\cup\{\aleph_0\})$, and

(4) $(\forall k \in \omega)(\exists m\in\omega)(\forall 
j\geq m)\allowbreak(k<f(j)$ and $k<g(j))$, and

(5) $(\forall\eta\in T)(\forall i\in{\rm dom}(\eta))(\eta(i)<f(i))$, and

(6) $(\forall n\in\omega)(\vert\{\eta(n)\,\colon\eta\in T$ and $n\in{\rm dom}(\eta)\}\vert\leq g(n))$.

\proclaim Definition 6.2. Suppose that $f$ and $g$ are functions as in Definition 6.1.
We say that $P$ is $(f,g)$-bounding iff\/ {\rm ${\bf 1}\forces_P``(\forall h\in {}^\omega\omega)\allowbreak [(\forall n\in\omega)\allowbreak(h(n)<f(n))$ implies $(\exists T\in V)\allowbreak
(T $ is an $(f,g)$-corseted tree and $h\in[T])]$.''}

\proclaim Lemma 6.3. Suppose $P$ is $(f^{g^k},g^{1/k})$-bounding for infinitely many $k\in\omega$,
and suppose $x<z$ are positive rational numbers. Suppose  $\gamma\in\omega$ and\/ {\rm ${\bf 1}\forces``T $ is an 
$(f^{g^\gamma}, g^x)$-corseted tree.''} Then ${\bf 1}\forces``(\exists H\in V)\allowbreak
(H$ is an $(f^{g^\gamma}, g^z)$-corseted tree and $T\subseteq H)$.''

Proof. Fix an integer $k$ such that $k> x$  and $P$ is $(f^{g^{\gamma+k}},g^{1/{(\gamma+k)}})$-bounding and $k>1/(z-x)$.
Let $X=\{n\in\omega\,\colon\allowbreak g(n)=\aleph_0\}$.

For every $m\in\omega-X$ define 

\medskip

\centerline{${\cal T}_m=\{S\subseteq\omega\,\colon\allowbreak
{\rm sup}(S)\leq f(m)^{g(m)^\gamma}$ and $\vert S\vert<g(m)^x\}$.}

\medskip

For every $m\in\omega-X$ define

\medskip

\centerline{${\cal T}'_m=\{i\in\omega\,\colon\allowbreak i\leq f(m)^{g(m)^{\gamma+k}}\}$.}

\medskip

Because $x<k$ we may choose, for each integer $m$ not in $X$, a one-to-one mapping $h_m$
 from ${\cal T}_m$
into ${\cal T}'_m$.

Define

\medskip

\centerline{${\cal T}=\{\xi\in{}^{<\omega}\omega\,\colon(\forall m\in \omega-X)(\xi(m)\in{\cal T}'_m)$}

\centerline{and $(\forall m\in X)(\xi(m)=1)\}$.}

\medskip

In $V[G_P]$ let $\zeta\in[{\cal T}]$ denote the function defined by

\medskip

\centerline{$(\forall m\in\omega-X)\allowbreak(\zeta(m)=h_m(\{\eta(m)\,\colon\allowbreak\eta\in T$ and $m\in{\rm dom}(\eta)\}))$}

\centerline{and $(\forall m\in X)\allowbreak(\zeta(m)=1)$.}

\medskip

Because $P$ is
$(f^{g^{\gamma+k}},g^{1/{(\gamma+k)}})$-bounding, 
we may take $H'\in V$ such that $H'$ is an $(f^{g^{\gamma+k}},g^{1/{(\gamma+k)}})$-corseted tree and
$\zeta\in[H']$. Define $H^*$ by

\medskip

\centerline{$H^*(m)=\bigcup\{h^{-1}_m(t)\,\colon\allowbreak(\exists\eta\in H')( t=\eta(m))$ and $t\in {\rm range}(h_m)\}$ for $m\in\omega-X$,}

\centerline{and $H^*(m)=\omega$ for $m\in X$.}

\medskip

When $g(m)$ is finite, we have

\medskip

{\centerline{$\vert H^*(m)\vert\leq\vert H'(m)\vert\cdot{\rm max}\{\vert h_m^{-1}(t)\,\colon\allowbreak t\in{\rm range}(h_m)\vert\}$}}

{\centerline{$\leq
 g^x(m)\cdot g^{1/{(\gamma+k)}}(m)< g^z(m)$.}}

\medskip

Let $H=\{\eta\in{}^{<\omega}\omega\,\colon(\forall i\in{\rm dom}(\eta))(\eta(i)\in H^*(i))\}$.
We have that $H$ is an $(f^{g^\gamma},g^z)$-corseted tree and ${\bf 1}\forces``T\subseteq H$.''  The Lemma is established.

\proclaim Lemmma 6.4. Suppose $n^*\in\omega$. Suppose $P$ is a forcing such that\/ {\rm
$V[G_P]\models``$for every countable $X\subseteq V$ there is a countable $Y\in V$ such that
$X\subseteq Y$.''} Suppose\/ {\rm $V[G_P]\models``\langle r_n\,\colon\allowbreak n\in\omega\rangle$ is a
bounded sequence of positive rational numbers and $y\in{\bf Q}$ and ${\rm sup}\{r_n\,\colon\allowbreak
n\in\omega\}<y$ and $(\forall n\in\omega)\allowbreak(T_n\in V$ is an 
$(f,g^{r_n})$-corseted tree).''}  
Then in  $V[G_P]$ there is an $(f,g^{y})$-corseted tree $T^*\in V$ and an increasing sequence of integers 
$\langle k_n\,\colon\allowbreak n\in\omega\rangle$ such that $k_0=0$ and $k_1\geq n^*$ and
$(\forall i>0)\allowbreak(i<k_i)$ and
for every $\eta\in{}^{<\omega}\omega$ we have 

\medskip

\centerline{$(\forall t\in{\rm dom}(\eta))(\exists j\in\omega)(\exists\nu\in T_{k_j})(k_j\leq t$ and $\nu(t)=\eta(t))$}

\centerline{implies $\eta\in T^*.$}

\medskip

Proof: The proof is similar to the proof of Lemma 4.7.
We note the following modifications. We must choose $x\in{\bf Q}$ such that ${\rm sup}\{r_n\,\colon\allowbreak n\in\omega\}<x<y$. 
By recursion choose $\langle k_n\,\colon\allowbreak n\in\omega\rangle$ an increasing sequence of integers such that $k_0=0$ and
$k_1\geq n^*$ and
$(\forall n>0)\allowbreak(\exists j\in\omega)\allowbreak(k_n\leq j$ implies
$(n+1)g(j)^{x}\leq g(j)^{y})$.

The definition of $T^*$ is changed to $T^*=\{\eta\in{}^{<\omega}\omega\,\colon(\forall t\in{\rm dom}(\eta))(\exists j\in\omega)(\exists\nu\in S'_{k_j})(k_j\leq t$ and $\nu(t)=\eta(t))\}$.

Clearly $T^*$ is a tree.

Claim. $T^*$ is an $(f,g^y)$-corseted tree.

Proof: Fix $t\in\omega$. 

Case 1:  $t\geq k_1$.

Choose $m\in\omega$ such that $k_m\leq t<k_{m+1}$.  We have that

\medskip

\centerline{
$\vert \{\eta(t)\,\colon\eta\in T^*$ and $t\in{\rm dom}(\eta)\}\vert=
\Sigma_{j\leq m}\vert \{\eta(t)\,\colon\allowbreak\eta\in T_{k_j}$ and
$t\in{\rm dom}(\eta)\}\vert$}

\centerline{$\leq (m+1)g^{x}(t)\leq g^{y}(t)$.}

\medskip

Case 2: $t<k_1$.

We have $\{\eta(t)\,\colon \eta\in T^*\}=\{\eta(t)\,\colon \eta\in T\}$, so it follows that
$\vert H(t)\vert\leq g^y(t)$.

The Claim is established.

The other requirements of the Lemma are the same as in the proof of Lemma 4.7. The Lemma is established.

\proclaim Lemma 6.5. Suppose $x<z$ are positive rational numbers, and suppose\/ {\rm
$(\forall n\in\omega)\allowbreak(T_n$ is an $(f,g^{x})$-corseted tree).}  Suppose $T$ is an $(f,g^x)$-corseted tree.
Then there is an $(f,g^{z})$-corseted tree $T^*\supseteq T$ and an increasing sequence of integers
$\langle m_i\,\colon\allowbreak i\in\omega\rangle$ such that 
for all $\eta\in T$ and all $i\in\omega$ and all $\nu\in T_{m_i}$ extending $\eta$, if
 ${\rm lh}(\eta)\geq m_i$, then
$\nu\in T^*$.

Proof. Choose a rational number $y$ such that $x<y<z$.  Choose $n^*\in\omega$ such that
$(\forall n\geq n^*)\allowbreak(2g^x(n)\leq g^y(n))$.

For every $n\geq n^*$ define
$T'_n=\{\eta\in{}^{<\omega}\omega\,\colon\allowbreak
\eta\in T$  or $\eta\in T_n$ and $\eta\restr n\in T\}$.
For every $n<n^*$ set $T'_n=T$.

It is easy to see that for every $n$ we have that $T'_n$ is an $(f,g^y)$-corseted tree.

By Lemma 6.4, we may take $T^*$ an $(f,g^{z})$-corseted tree, and an increasing sequence of integers
$\langle k_i\,\colon\allowbreak i\in\omega\rangle$, such that $k_0=0$ and $k_1\geq n^*$ and $(\forall i>0)\allowbreak(k_i>i)$ and
 for every $\eta\in{}^{<\omega}\omega$ we have

\medskip

\centerline{
$(\forall t\in{\rm dom}(\eta))\allowbreak(\exists j\in\omega)\allowbreak(\exists\zeta\in T'_{k_j})(k_j\leq t$ and $\zeta(t)=\eta(t))$}

\centerline{
  implies $\eta\in T^*$.}

\medskip

It is clear that $T\subseteq T^*$.

For every $i\in\omega$ set $m_i=k_{i+1}$.

Now suppose $\eta\in T$ and $i\in\omega$ and ${\rm lh}(\eta)\geq m_i$ and $\nu\in T_{m_i}$ and
$\nu$ extends $\eta$.  Because $\nu$ extends an element of $T$ of length at least $m_i$, we have
that $\nu\in T'_{m_i}$.
Choose $h\in[T'_{m_i}]$ such that $\nu$ is an initial segment of $h$.
It suffices to show that $h\in[T^*]$. 
Therefore it suffices to show that for every $t\in\omega$ that the following holds:

\medskip

$(*)_t$ $(\exists j\in\omega)\allowbreak(\exists\zeta\in T'_{k_j})(k_j\leq t$ and $\zeta(t)=h(t))$

\medskip

 Fix $t\in\omega$.

Fix $m^*\in\omega$ such that
$k_{m^*}\leq t<k_{{m^*}+1}$.

Case 1: $i<m^*$.

To see that $(*)_t$ holds, set $j=m^*$ and $\zeta=h\restr k_{m^*}$.
 We have $h\in[T'_{k_{i+1}}]\subseteq[T'_{k_{m^*}}]$,
so $\zeta\in T'_{k_{j}}$ and $k_{j}\leq t$. This completes Case 1.

Case 2: $m^*\leq i$.

We have $\eta\in T=T'_0$. Set $\zeta=\eta$ and $j=0$. Because $h$ extends $\eta$ we have $\zeta(t)=h(t)$, and so $(*)_t$ holds.
This completes Case 2.

The Lemma is established.

\proclaim Lemma 6.6. Suppose
$y$ and $z$ are positive rational numbers, and suppose
$\langle x_n\,\colon\allowbreak n\in\omega\rangle$ is a sequence of positive rational numbers
 such that $(\forall n\in\omega)\allowbreak
(x_n< x_{n+1}< y< z)$.  Suppose $T$ is an $(f,g^{x_0})$-corseted tree.
Suppose for every $n\in\omega$, we have $x_n^*\in{\bf Q}$ and 
$x_n^*< x_n< x^*_{n+1}$, and for each $n\in\omega$ we have
$\langle x_{n,j}\,\colon\allowbreak n\in\omega\rangle$ is a
 sequence of rational numbers such that
for every  $j\in\omega$ we have
$x_n< x_{n,j}<x_{n,j+1}< x_{n+1}^*$.
Suppose  $\langle T_{n,j}\,\colon\allowbreak n\in\omega$, $j\in\omega\rangle$ is
a sequence such that for every $n\in\omega$ and
 $j\in\omega$ we have that $T_{n,j}$ is an $(f,g^{x_{n,j}})$-corseted tree.
  Then there are $\langle T^n\,\colon\allowbreak n\in\omega\rangle$ and
$T^*$ such that $T^*$ is an $(f,g^{z})$-corseted tree and for every forcing notion $P$ we have  in $V[G_P]$ that 
$T\subseteq T^*$ and for every
$n\in\omega$ we have

(i) $T^n\subseteq T^{n+1}$ and $T^n$ is an
 $(f,g^{x_n})$-corseted tree, and

(ii)  for every $j\in\omega$ and every $g\in[T_{n,j}]\cap V[G_P]$ 
there is $k\in\omega$ such that for every $\eta\in T_{n,j}$ extending $g\restr k$, if 
 $\eta\restr k \in T^n\cap T^*$ then
$\eta\in T^{n+1}\cap T^*$.

\medskip

Proof:    Let $T^0=T$. 
Given $T^n$, build $\langle T'_{n,j}\,\colon\allowbreak j\in\omega\rangle$
as follows.  Let $T'_{n,0}=T^n$.  Given
$T'_{n,j}$ choose $m(n,j)\in\omega$ such that

\medskip

\centerline{$(\forall t\geq m(n,j))\allowbreak(2g^{x_{n,j}}(t)\leq g^{x_{n,j+1}}(t))$.}

\medskip

Set 

\medskip

\centerline{$T'_{n,j+1}=T'_{n,j}\cup\{\eta\in T_{n,j}\,\colon\allowbreak\eta\restr m(n,j)\in
T'_{n,j}\}$.}

\medskip

Claim 1.  Whenever $i\leq j<\omega$ we have $T'_{n,i}\subseteq T'_{n,j}$.

Proof.  Clear.

Claim 2. Suppose $T^n$ is an $(f,g^{x_n})$-corseted tree.  
Then $(\forall j\in\omega)\allowbreak(T'_{n,j}$
is an $(f,g^{x_{n,j}})$-corseted tree).

Proof:  It is clear that $T'_{n,0}$ is an $(f,g^{x_{n,0}})$-corseted tree.  
Assume that $T'_{n,j}$ is an $(f,g^{x_{n,j}})$-corseted tree.  Fix $t\in\omega$.

Case 1: $t<m(n,j)$.

We have that $\{\eta(t)\,\colon\allowbreak t\in T'_{n,j+1}$ and $t\in{\rm dom}(\eta)\}=
\{\eta(t)\,\colon\allowbreak\eta\in T'_{n,j}$ and $t\in{\rm dom}(\eta)\}$ and so 

\medskip

\centerline{$\vert \{\eta(t)\,\colon\allowbreak t\in T'_{n,j+1}$ and $t\in{\rm dom}(\eta)\}\vert\leq g^{x_{n,j}}(t)\leq
g^{x_{n,j+1}}(t)$.}

\medskip

Case 2: $t\geq m(n,j)$.  

We have $\{\eta(t)\,\colon\allowbreak t\in T'_{n,j+1}$ and $t\in{\rm dom}(\eta)\}
\subseteq(\{\eta(t)\,\colon\allowbreak t\in T'_{n,j}$ and $t\in{\rm dom}(\eta)\}\cup
\{\eta(t)\,\colon\allowbreak t\in T_{n,j}$ and $t\in{\rm dom}(\eta)\} )$.  Therefore
we have 

\medskip

\centerline{$\vert \{\eta(t)\,\colon\allowbreak t\in T'_{n,j+1}$ and $t\in{\rm dom}(\eta)\}\vert\leq
 2g^{x_{n,j}}(t)\leq g^{x_{n,j+1}}(t)$.}

\medskip

  The
Claim is established.

For each $n\in\omega$, using Claim 2 and Lemma 6.4 we my find an increasing sequence of integers
 $\langle k_{n,j}\,\colon\allowbreak j\in\omega\rangle$ and $T^{n+1}$ such that
$k_{n,0}=0$ and  $(\forall j>0)\allowbreak(k_{n,j}>j)$ and
if $T^n$ is an $(f,g^{x_n})$-corseted tree, then $T^{n+1}$ is
an $(f,g^{x_{n+1}})$-corseted tree such that
for all $\eta\in{}^{<\omega}\omega$, we have 

\medskip

\centerline{$(\forall t\in{\rm dom}(\eta))(\exists j\in\omega)(\exists\eta'\in T'_{n,k_{n,j}})(k_{n,j}\leq t$ and $\eta(t)=\eta'(t))$}

{\centerline{implies $\eta\in T^{n+1}$.}}

\medskip

This completes the construction of $\langle T^n\,\colon n\in\omega\rangle$ and $\langle T_{n,j}\,\colon
\allowbreak j\in\omega$, $n\in\omega\rangle$.

Applying mathematical induction, we have that each $T^n$ is in fact an $(f,g^{x_n})$-corseted tree.

Claim 3. $T^n\subseteq T^{n+1}$ for every $n\in\omega$.

Proof: By Claim 1 we have that $T^n\subseteq T'_{n,i}$ for every $i\in\omega$.
By the definition of $T^{n+1}$ we have that
$T^{n+1}\supseteq\bigcap\{T'_{n,k_{n,i}}\,\colon\allowbreak i\in\omega\}\supseteq
\bigcap\{T'_{n,i}\,\colon\allowbreak i\in\omega\}\supseteq T^n$. The Claim is established.

Applying Lemma 6.4 again we obtain an increasing sequence of integers $\langle k_n\,\colon\allowbreak n\in\omega\rangle$ and a
$(f,g^z)$-corseted tree $T^*$ such that $k_0=0$ and $(\forall n>0)\allowbreak
(n<k_n)$ and for every $\eta\in{}^{<\omega}\omega$,
we have that

\medskip

{\centerline{$(\forall t\in{\rm dom}(\eta))(\exists j\in\omega)(\exists\eta'\in T^{k_j})(k_j\leq t$ and $\eta(t)=\eta'(t))$}

{\centerline{implies $\eta\in T^*$.}}

\medskip

Notice  that $T^0\subseteq\bigcap\{T^n\,\colon\allowbreak n\in\omega\}\subseteq T^*$.

Now we verify that $\langle T^n\,\colon\allowbreak n\in\omega\rangle$ and $T^*$ satisfy the
remaining conclusions of the Lemma.  Accordingly, fix $P$ a forcing notion and work in $VG_P]$. Fix
$n\in\omega$ and $j\in\omega$ and
$g\in[T_{n,j}]$.
Let 

\medskip

\centerline{$k={\rm max}(k_{n+1},\allowbreak{\rm max}\{k_{n,j'}\,\colon\allowbreak j'\leq j+1\},
\allowbreak{\rm max}\{m(n,j')\,\colon\allowbreak
j'\leq j+1\})$.}

\medskip

  Fix $\eta\in T_{n,j}$ extending $g\restr k$ and assume that
 $\eta\restr k\in T^n\cap T^*$.

Claim 4. $\eta\in T^{n+1}$.  

Fix $t\in{\rm dom}(\eta)$.  

Let $j'$ be the unique integer such that $k_{n,j'}\leq t<k_{n,j'+1}$.

It suffices to show 

\medskip

$(*)_t$  $(\exists j^*\leq j')(\exists\eta'\in T'_{n,k_{n,j^*}})(\eta(t)=\eta'(t))$.

\medskip

Case 1:  $j> j'$.

Let $\eta'=\eta\restr k_{n,j'+1}$ and let $j^*=j'$.
Because $\eta'$ is an initial segment of
$\eta\restr k\in T^n\subseteq T'_{n,k_{n,j'}}$, we have that
  $(*)_t$ holds.

Case 2: $j\leq j'$.

Let $\eta'=\eta$ and $j^*=j$.
Because $m(n,j)\leq k$ and $g\restr k=\eta\restr k\in T^n$, we have $g\restr m(n,j)\in T^n\subseteq T'_{n,j}$.
Because we also have $\eta\in T_{n,j}$, we conclude  that $\eta'\in T'_{n,j+1}\subseteq T'_{n,k_{n,j}}$.
It is easy to see that $(*)_t$ holds.

Claim 4 is established.

Claim 5. $\eta\in T^*$.

Fix $t\in{\rm dom}(\eta)$.  

Let $i$ be the unique integer such that $k_{i}\leq t<k_{i+1}$.

It suffices to show 

\medskip

\centerline{  $(\exists i'\leq i)(\exists\eta'\in T^{k_{i'}})(\eta(t)=\eta'(t))$.}

\medskip

Case 1: $i<n$.

Because $\eta\restr k\in T^*$ we have $g\restr k_{i+1}\in T^*$ and so we may take
$i'\leq i$ and $\eta'\in T^{k_{i'}}$ such that $g(t)=\eta'(t)$.

Case 2: $n\leq i$.

Let $i' = i$ and $\eta'=\eta\restr k_{i+1}$. We have
$\eta'\in T^n\subseteq T^{k_i}$.

Claim 5 and the Lemma are established.

The proof of the following Theorem is very similar to the proof of Lemma 4.10.

\proclaim Theorem 6.7.  Suppose $\langle P_\eta\,\colon\allowbreak\eta\leq\kappa\rangle$ is a 
countable support iteration based on $\langle Q_\eta\,\colon\allowbreak\eta<\kappa\rangle$ and suppose\/ {\rm $(\forall\eta<\kappa)\allowbreak({\bf 1}\forces_{P_\eta}``$for infinitely
many $k\in\omega$ we have that $Q_\eta$ 
is proper and $(f^{g^k},g^{1/k})$-bounding''),} and suppose $\gamma\in\omega$.
Suppose $\lambda$ is a sufficiently large regular cardinal 
and $\alpha<\kappa$ 
and  $x$ and $z$ are $P_\alpha$-names and $T$ is a $P_\kappa$-name  
and\/ {\rm ${\bf 1}\forces_{P_\alpha}``x$ and $z$ are positive rational numbers and $x<z$''
  and ${\bf 1}\forces_{P_\kappa}``T$ is an $(f^{g^\gamma}, g^x)$-corseted tree.''}
Suppose $N$ is a countable elementary submodel of $H_\lambda$ and $\{P_\kappa, \alpha, f,g,  x, z, T
\}\in N$. 
Suppose
 $p\in P_\alpha$ and $p$ is $N$-generic.
Then\/ {\rm
$p\forces``(\forall q\in P_{\alpha,\kappa}\cap N[G_{P_\alpha}])\allowbreak(\exists q'\leq q)\allowbreak
(\exists H)\allowbreak (H$ is an $(f^{g^\gamma},g^z)$-corseted tree and
$ q'
\forces`T\subseteq H$').''}

Proof: The proof proceeds by induction on $\kappa$.  We assume that 
  $\lambda$, 
$N$, $\alpha$, $p$,  $x$, $z$, and $T$ are as in the hypothesis of the Theorem. Let $f'=f^{g^\gamma}$.

Fix $q$ a $P_\alpha$-name in $N$ such that ${\bf 1}\forces``q\in P_{\alpha,\kappa}$.''

Case 1. $\kappa=\beta+1$.

Fix $y$ a $P_\alpha$-name in $N$ such that ${\bf 1}\forces``y$ is rational and $x< y< z$.''

Using Lemma 6.4 choose $\tilde q$ and $H'$ such that

\medskip

\centerline{ 
${\bf 1}\forces_{P_\beta}``\tilde q\leq q(\beta)$ and $H'$ is an $(f',g^y)$-corseted tree and
$\tilde q\forces``T\subseteq H'$.'\thinspace''}

\medskip

We may assume that the names $\tilde q$ and $H'$ are elements of $N$.
Use the induction hypothesis to choose  $P_\alpha$-names $q^*$ and
 $H$  such that

\medskip

\centerline{
$p\forces``q^*\leq q\restr\beta$ and $H$ is an $(f',g^x)$-corseted tree and
$q^*\forces`H'\subseteq H$.'\thinspace''}

\medskip

We have that $p\forces``(q^*,\tilde q)\forces`T\subseteq H$.'\thinspace''  Case 1 is established.

\medskip

Case 2.  ${\rm cf}(\kappa)>\omega$.

Because no $\omega$-sequences of ordinals can be added at limit stages of uncountable cofinality,
 we may take $\beta$ and $T'$ and $q'$ to be $P_\alpha$-names in $N$ such that

\medskip

\centerline{
${\bf 1}\forces``\alpha\leq\beta<\kappa$ and $T'$ is a
$P_{\alpha,\beta}$-name and $q'\leq q$ and}

\centerline{
${\bf 1}\forces_{P_{\alpha,\beta}}`T'$ is an $(f',g^x)$-corseted tree' and
$q'\forces_{P_{\alpha,\kappa}}`T'=T$.'\thinspace''}

\medskip

For every $\beta_0\in\kappa\cap N$ such that $\alpha\leq\beta_0$ 
let $\tilde q(\beta_0)$ and $H(\beta_0)$  be  $P_\alpha$-names in $N$ such that

\medskip

\centerline{
${\bf 1}\forces``$if $\beta=\beta_0$ and there is some $\tilde q\leq q'\restr\beta$ and some 
$H^*$}

\centerline{such that $H^*$ is an $(f',g^z)$-corseted tree and
$\tilde q\forces`T'\subseteq H^*$,'}

\centerline{
then $q^*(\beta_0)$ and $H(\beta_0)$  are witnesses thereto.''}

\medskip

Let $ q^*$ and $H$  and $s$ be $P_\alpha$-names such that for every $\beta_0
\in \kappa\cap N$, if $\alpha\leq\beta_0$, then

\medskip

\centerline{
${\bf 1}\forces``$if $\beta=\beta_0$ then $ q^*=q^*(\beta_0)$ and $H=H(\beta_0)$
 and 
$s\in P_{\alpha,\kappa}$}

\centerline{and $s\restr\beta=q^*$ and
$s\restr[\beta,\kappa)=q'\restr[\beta,\kappa)$.''}

\medskip

Claim 1: $p\forces``s\leq q$ and $s\in N[G_{P_\alpha}]$ and $s
\forces`T\subseteq H$.'\thinspace''

Proof:  Suppose $p'\leq p$. Take $p^*\leq p'$ and $\beta_0<\kappa$ such that
$p^*\forces``\beta_0=\beta$.''  Because the name $\beta$ is in $N$ and $p^*$ is $N$-generic,
we have that $\beta_0\in N$.  Notice by the induction hypothesis
we have

\medskip

\centerline{$p\forces``$there is some $q^\#\leq q'\restr\beta_0$ and 
some $(f' g^z)$-corseted tree
$H^\#$ }

\centerline{
such that $q^\#\forces`T'\subseteq H^\#$.'\thinspace''}

\medskip

\noindent Hence 

\medskip

\centerline{$p^*\forces``q^*=q^*(\beta_0)\leq q'\restr\beta$ and $H=H(\beta_0)$ and
$H$ is an $(f',g^z)$-corseted}

\centerline{tree and $q^*\forces`T'\subseteq H$ and
$q'\restr[\beta,\kappa)\forces``T'=T$.''\thinspace'\thinspace''}

\medskip

\noindent  Therefore
$p^*\forces``s\forces`T\subseteq H$.'\thinspace''

Claim 1 is established.
This completes Case 2.

Case 3. ${\rm cf}(\kappa)=\omega$.

Let $\langle\alpha_n\,\colon\allowbreak n\in\omega\rangle$ be an increasing 
sequence from $\kappa\cap N$ cofinal in $\kappa$ such that 
$\alpha_0=\alpha$.

In $V[G_P]$, let $X$, $\langle h_m\,\colon\allowbreak m\in\omega-X\rangle$, $\zeta$,
 $\langle {\cal T}_m\,\colon\allowbreak m\in\omega_X\rangle$, $\langle {\cal T}'_m\,\colon\allowbreak m\in\omega_X\rangle$,
${\cal T}$, and $H'$ be as in the proof of Lemma
6.3 with $P_{\alpha,\kappa}$ playing the role of $P$.  
We may assume each of these $P_\alpha$-names are in $N$.

Using Lemma 3.1, fix
 $\langle (p_n,\zeta_n)\,\colon\allowbreak n\in\omega\rangle\in N$ 
(that is, the sequence of names is an element of $N$ but not necessarily their values) such that 
${\bf 1}\forces``p_0\leq q$'' and 
 for every $n\in\omega$ we have that each of the following holds:

(0) $p_n$ is a $P_\alpha$-name for an element of $P_{\alpha,\kappa}$, and

(1) For every $k\leq n$ we have ${\bf 1}\forces_{P_{\alpha_n}}``p_0\restr[\alpha_n,\kappa)
\forces`(\forall t<k)(\zeta(t)=\zeta_n(t))$,'\thinspace'' and

(2) $\zeta_n$ is a $P_{\alpha_n}$-name for an element of $[{\cal T}]$, and

(3) ${\bf 1}\forces_{P_{\alpha_n}}``p_0\restr[\alpha_n,\alpha_{n+1})\forces`(\forall t<k)(\zeta_n(t)
=\zeta_{n+1}(t))$ for every $k\leq n+1$,'\thinspace'' and

(4) ${\bf 1}\forces_{P_\alpha}``p_{n+1}\leq p_n$,'' and

(5) whenever $k\leq m<\omega$ we have 

\medskip

\centerline{${\bf 1}\forces_{P_{\alpha_n}}``p_m\restr[
\alpha_n,\alpha_{n+1})\forces`(\forall t<k)(\zeta_n(t)=\zeta_{n+1}(t))$.'\thinspace''}

\medskip

Claim 2.  Suppose $\alpha\leq\beta\leq\zeta<\kappa$ and ${\bf 1}\forces_{P_\alpha}
``x'< z'$ are positive rational numbers''  and
suppose $T'$ is a $P_\zeta$-name for an $(f',g^{x'})$-corseted tree. Then

\medskip

\centerline{${\bf 1}\forces_{P_\beta}``V[G_{P_\alpha}]\models`(\forall q\in P_{\beta,\zeta})(\exists q^*\leq q')(\exists H)$}

\centerline{$( 
H$ is an $(f',g^{z'})$-corseted tree and
$q^*\forces``T'\subseteq H$'').'\thinspace''}

\medskip

Proof: Given $r_1\in P_\alpha$ and a $P_\alpha$-name $r_2$ for an element of $P_{\alpha,\beta}$ and a 
$P_\beta$-name $q$ for an element of $P_{\beta,\zeta}\cap V[G_{P_\alpha}]$, choose $\lambda'$ a sufficiently large
regular cardinal and $N'$ a countable elementary substructure of $H_{\lambda'}$ containing $\{r_1,r_2,q,P_\kappa,\alpha,\beta,\zeta,x',z', T'\}$.
Choose $r_1'\leq r_1$ such that $r_1'$ is $N'$-generic.  By the overall induction hypothesis (i.e., because $\zeta<\kappa$)
we have

\medskip

\centerline{$r'_1\forces``(\exists s\leq(r_2,q))(\exists H)(H$ is an $(f',g^{z'})$-corseted tree and $s\forces`T'\subseteq H$').''}

\medskip

Consequently, we may choose $s$ and $h$ such that

\medskip

\centerline{$(r'_1,s\restr\beta)\forces``V[G_{P_\alpha}]\models`s\restr[\beta,\zeta)\leq q$ and}

\centerline{$H$ is an $(f',g^{z'})$-corseted
 tree and $s\restr[\beta,\zeta)\forces`T'\subseteq H$').''}

\medskip

The Claim is established.

In $V[G_{P_\alpha}]$ fix a positive rational number $y< 1/(\gamma+k)$ where $k$ is as in the proof of Lemma 6.3.
  We may assume that the name $y$ is in $N$.

Let $\Omega=\{x'\in N\,\colon\allowbreak x'$ is a $P_\alpha$-name and 
${\bf 1}\forces``x'$ is rational and $0< x'< y$''$\}$. Let $\langle y_n\,\colon\allowbreak n\in\omega\rangle$ enumerate $\Omega$.  Build $\langle x^*_n\,\colon\allowbreak n\in\omega\rangle$ as follows.  Let $x^*_0=y_0$, and for each $n\in\omega$ 
choose $x^*_{n+1}\in \Omega$ such that
${\bf 1}\forces``x^*_n< x^*_{n+1}$ and $y_{n+1}< x^*_{n+1}$.''
Also build $\langle x_n\,\colon\allowbreak n\in\omega\rangle$ a sequence of 
elements of $\Omega$ such that for each $n\in\omega$ we have
$x^*_n< x_n< x^*_{n+1}$.

For each $n\in\omega$ let $\langle T_{n,j}\,\colon\allowbreak j\in\omega\rangle$ list all $P_{\alpha}$-names $T'\in N$ such that
in $V[G_{P_\alpha}]$ we have for some $y'< x^*_{n+1}$ we have that $T'\subseteq{\cal T}$ is an
$(f',g^{y'})$-corseted tree, and build $\langle x_{n,j}\,\colon\allowbreak j\in\omega\rangle$ a sequence of elements of $\Omega$
 such that in $V[G_{P_\alpha}]$ for every $j\in\omega$ we have that
$x_n< x_{n,j}< x_{n,j+1}< x^*_{n+1}$ and $T_{n,j}\subseteq{\cal T}$ is an $(f',g^{x_{n,j}})$-corseted tree.

 Using Lemma 6.6, choose $T^*\subseteq {\cal T}$ an $(f',g^y)$-corseted tree and
 $\langle T^n\,\colon\allowbreak n\in\omega\rangle\in V[G_{P_\alpha}]$
 a sequence such that $T^*\in V[G_{P_\alpha}]$ and $T^0\subseteq T^*$ and
$\zeta\in[T^0]$ and for every $n\in\omega$ we have 
$T^n\subseteq{\cal T}$ is an $(f',g^{x_n})$-corseted tree and $T^n\subseteq T^{n+1}$ and, in $V[G_{P_\kappa}]$, we have that
 for every $j\in\omega$ and every 
$g\in[T_{n,j}]$  there is $k\in\omega$ such that for
every $\eta\in T_{n,j}$ extending $g\restr k$, if
 $\eta\restr k\in T^n\cap T^*$ then
$\eta\in T^{n+1}\cap T^*$.  

Note that the reason we worked in $V[G_{P_\kappa}]$ in the preceding paragraph is because we wish to allow $g$ to
range over $[T_{n,j}]$ with the brackets interpreted in $V[G_{P_\kappa}]$.

We may assume that the names $T^*$ and $\langle T^n\,\colon\allowbreak n\in\omega\rangle$ are in $N$.

For each $n\in\omega$ let $F_{n,0}$ and $F_{n,2}$ and $y^*_n$ be $P_{\alpha_n}$-names  such that

\medskip

\centerline{
 ${\bf 1}\forces``F_{n,0}$ and
$F_{n,2}$ and $y^*_n$ are functions, all three of which are in $V[G_{P_\alpha}]$}

\centerline{each of whose domain is equal to $P_{\alpha_n,\alpha_{n+1}}$, such that}

\centerline{ $(\forall q'\in P_{\alpha_n,\alpha_{n+1}}\cap V[G_{P_\alpha}]
 )\allowbreak (F_{n,0}(q')\subseteq {\cal T}$ and}

\centerline{
$F_{n,0}(q')\in V[G_{P_\alpha}]$ is an $(f',g^{y^*_n(q')})$-corseted tree}

\centerline{
 and $y^*_n(q')$ is rational and  $0<y^*_n(q')< x^*_{n+1}$ and
$F_{n,2}(q')\leq q'$}

\centerline{ and $F_{n,2}(q')\forces`\zeta_{n+1}\in [F_{n,0}(q')]$'.''}

\medskip

 We may assume that the names $F_{n,0}$ and
$F_{n,2}$ and $y^*_n$ are in $N$. 

For each $n\in\omega$ we may, in $V[G_{P_{\alpha_n}}]$, use Lemma 6.6 to choose a positive rational number $y_n< x_{n+1}$
and
$\tilde T_n\subseteq{\cal T}$ an $(f',g^{y_n})$-corseted tree  and
$\langle k^n_i\,\colon\allowbreak i\in\omega\rangle\in V[G_{P_\alpha}] $ an increasing sequence of
integers 
such that $\tilde T_n\in V[G_{P_\alpha}]$ and
 $T^n\subseteq\tilde T_n$ and for every $\eta\in T^n$ and every $i\in\omega$ and every
$\nu\in F_{n,0}(p_{k^n_i}\restr[\alpha_n,\alpha_{n+1}))$,
 if 
${\rm length}(\eta)\geq k^n_i$ and 
$\nu$ extends $\eta$, then $\nu\in\tilde T_n$. 

\medskip

We may assume that for each $n\in\omega$ the $P_{\alpha_n}$-names $\tilde T_n$ and $\langle k^n_i\,\allowbreak i\in\omega\rangle$ are in $N$.

Claim 3.  We may be build $\langle r_n\,\colon\allowbreak n\in\omega\rangle$ such that
$r_0=p$  and for every $n\in\omega$ we have that the following hold:

(1) $r_n\in P_{\alpha_n}$ is $N$-generic, and

(2) $r_{n+1}\restr\alpha_n=r_n$, and

(3) $r_n\forces``\zeta_n\in [T^n]\cap [T^*]$,'' and

(4) $p\forces``r_n\restr[\alpha,\alpha_n)\leq p_0\restr\alpha_n$.''

The proof of this Claim is the same as the proof of Claim 3 of the proof of Lemma 4.10.  The Claim is established.

Let $q'$ be a $P_\alpha$-name such that

\medskip

\centerline{
$p\forces``q'=\bigcup\{r_n\restr[\alpha,\alpha_n)\,\colon\allowbreak n\in\omega\}$.''}

\medskip

 Define $H^*$ by
$H(m)=\bigcup\{h^{-1}_m(t)\,\colon\allowbreak t\in \{\eta(m)\,\colon \allowbreak\eta\in T^*$ and
$m\in{\rm dom}(\eta)\}\cap{\rm range}(h_m)\}$ for $m\in\omega-X$,
and $H^*(m)=\omega$ for $m\in X$.

Let $H=\{\eta\in{}^{<\omega}\omega\,\colon(\forall i\in{\rm dom}(\eta))(\eta(i)\in H^*(i))\}$.

As in Lemma 6.3, we have that $H$ is an $(f^{g^\gamma},g^z)$-corseted tree and ${\bf 1}\forces``T\subseteq H$.''  
By Claim~3, we have that $q'$ and $H$ satisfy the requirements of the Theorem.

The Theorem is established.

\proclaim Corollary 6.8.   Suppose $\langle P_\eta\,\colon\allowbreak\eta\leq\kappa\rangle$ is a countable support iteration
based on $\langle Q_\eta\,\colon\allowbreak\eta<\kappa\rangle$ and suppose that for every $\eta<\kappa$ we have that\/ {\rm
${\bf 1}\forces``$for infinitely many $k\in\omega$ we have that $Q_\eta$ is proper and
$(f^{g^k},g^{1/k})$-bounding.''} Then $P_\kappa$ is
$(f^{g^k},g^{1/k})$-bounding for every positive $k\in\omega$.

Proof: By Theorem 6.7 with $\alpha=0$.

\section{$P$-point property}

In this section we define the $P$-point property and prove that it is preserved by countable support iteration of proper
forcings.  This is due to Shelah [12, Conclusion VI.2.12G].

\proclaim Definition 7.1. Suppose $n\in\omega$ and $x\in{}^\omega(\omega-\{0\})$ is strictly increasing.
We say that $(j,k,m)$ is an $x$-bound system above $n$ iff each of the following holds:

(1)   $j$ and $m$ are functions from $k+1$ into $\omega$, and

(2) $j(0)>x(n+m(0)+1)$, and

(3) $(\forall l< k)\allowbreak(j(l+1)>x(j(l)+m(l+1)+1))$.

\proclaim Definition 7.2. Suppose $n\in\omega$ and $x\in{}^\omega(\omega-\{0\})$ is strictly increasing
and $(j,k,m)$ is an $x$-bound system above $n$ and $T$ is a tree.  We say that $T$ is a
$(j,k,m,\eta)$-squeezed tree iff each of the following holds:

(1) ${\rm dom}(\eta)=\{(l,t)\in\omega^2\,\colon\allowbreak
l\leq k$ and $t\leq m(l)\}$, and

(2) $(\forall (l,t)\in{\rm dom}(\eta))\allowbreak
(\eta(l,t)\in {}^{j(l)}\omega)$, and

(3) $(\forall\nu\in T)\allowbreak(\exists (l,t)\in{\rm dom}
(\eta))\allowbreak(\nu$ is comparable with $\eta(l,t))$.

It is easy to see that the following Definition is equivalent to [12, Definition VI.2.12A].

\proclaim Definition 7.3. We say that $T$ is $x$-squeezed iff for every $n\in\omega$ there is some $x$-bound system
$(j,m,k)$ above $n$ such that $T$ is $(j,k,m,\eta)$-squeezed for some $\eta$.

In other words, $T$ is $x$-squeezed when, living above any given level of $T$, say $\{\xi\in T\,\colon\allowbreak
{\rm lh}(\xi)=n+1\}$, there is a maximal antichain ${\cal A}$ of $T$ that can be decomposed
as ${\cal A}=\bigcup\{{\cal A}_l\,\colon\allowbreak l\leq k\}$ where each
${\cal A}_l$ is a subset of $\{\xi\in T\,\colon\allowbreak{\rm lh}(\xi)=j(l)\}$ of cardinality at most $m(l)+1$, such that
the levels of ${\cal A}$ are stratified so sparsely that conditions (2) and (3) of Definition 7.1 hold. 
 Notice that for any given $l\leq k$ we may have that
$\{\eta(l,t)\,\colon\allowbreak t\leq m(l)\}$ is a proper superset of ${\cal A}_l$; indeed, it need not even be a subset of $T$.
We could modify Definition 7.2 to require this, but there is no need to do so.

\proclaim Lemma 7.4. Suppose $1\ll x\ll y$ and both $x$ and $y$ are strictly increasing
and $T$ is a $y$-squeezed tree. Then $T$ is an $x$-squeezed tree.

Proof: Every $y$-bound system is an $x$-bound system.

\proclaim Definition 7.5.  We say that $P$ has the $P$-point property iff for every $x\in{}^\omega(\omega-\{0\})$ strictly increasing, we have

{\centerline{${\bf 1}\forces``(\forall f\in{}^\omega\omega)(\exists H\in V)(f\in [H]$ and $H$ is an $x$-squeezed tree).''}}

\medskip

\proclaim Lemma 7.6.  $P$ has the $P$-point property iff for every $x\in
{}^\omega(\omega-\{0\})$  strictly increasing and every
 $p\in P$, if\/ {\rm $p\forces_P``f\in{}^\omega\omega$''}  there are $q\leq p$ and an 
 $x$-squeezed tree
 $H$ such that 
$q\forces``f\in[H]$.''

Proof: Assume that $P$ has the $P$-point property.  Given $x$, $p$, and $f$, there is $q\leq p$ and $H\subseteq{}^{<\omega}\omega$ such that
$q\forces``f\in[H]$ and $H$ is an $x$-squeezed tree.'' By the Shoenfield Absoluteness Theorem we have that
$H$ is an $x$-squeezed tree.

The other direction is immediate, and so the Lemma is established.

\proclaim Lemma 7.7. Suppose $T$ is an $x$-squeezed tree and $n\in\omega$. Then $T\cap{}^n\omega$ is finite.

Proof. Fix  $(j,k,m)$ an $x$-bound system above $n$ and fix $\eta$ such that
$T$ is a $(j,k,m,\eta)$-squeezed tree.
We have $T\cap{}^n\omega\subseteq\{\eta(s,t)\restr n\,\colon\allowbreak
t\leq j(k)$ and $s\leq m(t)\}$.

\proclaim Lemma 7.8. Suppose that $P$ has the $P$-point
property. Then $P$ is ${}^\omega\omega$-bounding.

Proof: Suppose $p\in P$ and $p\forces``f\in{}^\omega\omega$.'' Pick $x\in{}^\omega
(\omega-\{0\})$ such that $1\ll x$, and take $q\leq p$ and $H$ an $x$-squeezed
tree such that 
$q\forces``f\in[H]$.''  By Lemma 7.7 we may define $h\in{}^\omega\omega$ by
$(\forall n\in\omega)\allowbreak(h(n)={\rm max}\{\nu(n)\,\colon\allowbreak
\nu\in H$ and $n\in{\rm dom}(\nu)\})$. Clearly $q\forces``f\leq h$,'' and
the Lemma is established.

\proclaim Lemma 7.9. Suppose $P$ has the Sacks property. 
Then $P$ has the $P$-point property.

Proof. Suppose $x\in{}^\omega(\omega-\{0\})$ is strictly increasing and $p\in P$ and
$p\forces``f\in{}^\omega\omega$.'' 
Choose $y\in{}^\omega(\omega-\{0\})$ monotonically non-decreasing such that for $n> x(3)$ 
we have that $y(n)$ is the greatest $t\in\omega$ such that $x(3t)<n$.
Using the Sacks property, choose $q\leq p$ and $H$ a $y$-sized tree such that $q\forces``f\in[H]$.''

Notice that for all $t>0$  we have $y(x(3t))$ is less than or equal to
 the greatest integer $k$
satisfying $x(3k)<x(3t)$, and therefore we have

\begin{itemize}

\item[$(*)$] $(\forall t\in\omega)(y(x(3t))<t)$.

\end{itemize}

Suppose $n>x(3)$.
Let $j$ be such that ${\rm dom}(j)=\{0\}$ and $j(0)=x(2n)+1$;
let $k=0$; and let $m$ be such that ${\rm dom}(m)=\{0\}$ and $m(0)=\vert H\cap
{}^{j(0)}\omega\vert$. 

Claim: $(j,k,m)$ is an $x$-bound system above $n$.

Proof: We have $x(n+m(0)+1)\leq x(n+1+y(x(2n)+1))\leq x(n+1+y(x(3n)))\leq x(n+1+n-1)<x(2n)+1=j(0)$. 
The first inequality is because $m(0)=\vert H\cap{}^{j(0)}\omega\vert\leq y(j(0))=y(x(2n)+1)$.
The second inequality is because $x$ is stricltly increasing and $y$ is monotonically non-decreasing.  The third
inequality is by $(*)$.

The Claim is established.

Define $\eta$ with domain equal to $\{(0.i)\,\colon\allowbreak
i<m(0)\}$ and such that $\langle \eta(0,i)\,\colon\allowbreak i<m(0)\rangle$
enumerates $H\cap{}^{j(0)}\omega$.  Clearly $H$ is a $(j,k,m,\eta)$-squeezed tree, so the Lemma is established.

\proclaim Lemma 7.10. Suppose $y\in{}^\omega(\omega-\{0\})$ is strictly increasing and $T$ and $T'$ are $y$-squeezed trees.
Then $T\cup T'$ is a $y$-squeezed tree.

Proof:  Given $n\in\omega$, choose $(j,k,m,\eta)$  such that
$(j,k,m)$ is a $y$-bound systems above $n$ and
$T$ is $(j,k,m,\eta)$-squeezed. Let $h=j(k)$. Choose $(j',m',k')$ a $y$-bound system above $h$ and choose $\eta'$
such that $T'$ is $(j',\allowbreak k' ,\allowbreak m',\eta')$-squeezed.
We proceed to fuse $(j,k,m,\eta)$ with $(j',\allowbreak k'\allowbreak m',\allowbreak
\eta')$. 
For every 
$l\leq k$ let $j^*(l)=j(l)$ and for every $l$ such that $k<l\leq k+k'+1$ let
$j^*(l)=j'(l-k_n-1)$.
Let $k^*=k+k'+1$.
For every 
$l\leq k$ let $m^*(l)=m(l)$ and for every $l$ such that $k<l\leq k+k'+1$ let
$m^*(l)=m'(l-k-1)$.
For every 
$l\leq k$ and $\beta\leq m(l)$ let $\eta^*(l,\beta)=\eta(l,\beta)$ and for every $l$ such that $k<l\leq k+k'+1$ and
every $\beta\leq m'(l-k-1)$ let
$\eta^*(l,\beta)=\eta'(l-k-1,\beta)$.
It is straightforward to verify that $(j^*,k^*,m^*)$ is a $y$-bound system above $n$ and
that $T\cup T'$ is $(j^*,k^*,m^*,\eta^*)$-squeezed.

The Lemma is established.

\proclaim Definition 7.11. Suppose  $n\in\omega$ and
$h\in{}^\omega\omega$ and  $y\in{}^\omega(\omega-\{0\})$ is strictly increasing.
Suppose $(j,m,k)$ is a $y$-bound system above $n$.  We say that $(j,m,k)$ is $h$-tight iff  $j(k)<h(n)$.
For $T$  a $y$-squeezed tree, we say that $T$ is $h$-tight iff for every $n\in\omega$ there is an
$h$-tight $y$-bound system
$(j,m,k)$  above $n$ such that 
 for some $\eta$ we have that $T$ is $(j,m,k,\eta)$-squeezed for some $\eta$.

\proclaim Lemma 7.12. Suppose $P$ has the $P$-point property and $y\in{}^\omega(\omega-\{0\})$ 
is strictly increasing and\/ {\rm ${\bf 1}\forces``T$ is a $y$-squeezed tree.''} 
Then\/ {\rm ${\bf 1}\forces``(\exists H\in V)\allowbreak(H$ is a $y$-squeezed tree and
$T\subseteq H$).''}

Proof. Suppose $p\in P$ and $p\forces``T$ is a $y$-squeezed tree.'' Fix $q'\leq p$.
 By Lemma 7.8 we may choose $h\in{}^\omega\omega$ and 
$q\leq q'$ such that 
$q\forces``T$ is $h$-tight.''

Define $z\in{}^\omega(\omega-\{0\})$ by $z(0)=0$ and

\medskip

{\centerline{$(\forall n\in\omega)\allowbreak(z(n+1)= h(z(n))$.}}

\medskip

For every $n\in\omega$ let ${\cal T}_n=\{t\subseteq{}^{< h(n)}\omega\,\colon\allowbreak 
t=T\cap{}^{< h(n)}\omega$ for some $h$-tight $y$-squeezed tree $T\}$. 

Let ${\cal T}=\bigcup\{{\cal T}_n\,\colon n\in\omega\}$. We implicitly fix an isomorphism
from ${}^{<\omega}\omega$ onto ${\cal T}$.

Using the fact that $P$ satisfies the $P$-point property, fix $q^*\leq q$ and ${\cal C}\subseteq{\cal T}$ such that ${\cal C}$ is a
$z$-squeezed tree and $q^*\forces``(\forall n \in\omega)\allowbreak(T\cap{}^{< h(n)}\omega\in{\cal C})$.''

Define $H^*=\bigcup{\cal C}$ and let $H=\{\nu\in H^*\,\colon\allowbreak(\forall n\in\omega)\allowbreak(\exists\eta\in{}^n\omega\cap H^*)
\allowbreak(\eta$ is comparable with $\nu)\}$.

Pick a $z$-bound system $(j^*,k^*,m^*)$ above $n$ and $\eta^*$ such that  ${\cal C}$ is a $(j^*,k^*,\allowbreak m^*,\allowbreak\eta^*)$-squeezed tree.

Fix $n\in\omega$. We show that there is a $y$-bound system $(j,m,k)$ above $n$ such that for some $\eta$ we have
that $H$ is $(j,m,k,\eta)$-squeezed.

Claim 1. For every $\beta\leq m^*(0)$ we have ${\rm ht}(\eta^*(0,\beta))\geq
h(z(n+\beta+1))$.  For every non-zero $\alpha\leq k^*$ and every $\beta\leq m^*(\alpha)$ we have ${\rm ht}(\eta^*(\alpha,\beta))\geq
h(z(j^*(\alpha-1)+\beta+1))$.

Proof:  For every $\beta\leq m^*(0)$ we have ${\rm ht}(\eta^*(0,\beta))=h({\rm rk}_{\cal T}(\eta^*(0,\beta)))=h(j^*(0))\geq
h(z(n+m^*(0)+1))\geq
h(z(n+\beta+1))$.  For every non-zero $\alpha\leq k^*$ and every 
$\beta\leq m^*(\alpha)$ we have ${\rm ht}(\eta^*(\alpha,\beta))=h({\rm rk}_{\cal T}(\eta^*(\alpha,\beta)))=h(j^*(\alpha))\geq
h(z(j^*(\alpha-1)+m^*(\alpha)+1))\geq
h(z(j^*(\alpha-1)+\beta+1))$.

By Claim 1 we may construct $y$-bound systems as follows.
For every  $\beta\leq m^*(0)$, fix an $h$-tight $y$-bound system $(j^{0,\beta},m^{0,\beta},\allowbreak
k^{0,\beta})$ above $z(n+\beta+1)$ along with $\eta^{0.\beta}$ such that
for some $(j^{0,\beta},m^{0,\beta},k^{0,\beta},\eta^{0,\beta})$-squeezed tree $T$ we have $\eta^*(0,\beta)=
{}^{<h(z(n+\beta+1))}\omega\cap T$, and for every non-zero $\alpha\leq k^*$ and $\beta\leq m^*(\alpha)$, 
fix an $h$-tight $y$-bound system $(j^{\alpha,\beta},m^{\alpha,\beta},\allowbreak
k^{\alpha,\beta})$ above $z(j^*(\alpha-1)+\beta+1)$ along with $\eta^{\alpha.\beta}$ such that
for some $(j^{\alpha,\beta},m^{\alpha,\beta},k^{\alpha,\beta},\eta^{\alpha,\beta})$-squeezed tree $T$ we have $\eta^*(\alpha,\beta)=
{}^{<h(z(j^*(\alpha-1)+\beta+1))}\omega\cap T$.

We define

\medskip

\centerline{$\hat\jmath(\alpha,\beta,\gamma)=j^{(\alpha,\beta)}(\gamma)$}

\medskip

\noindent and 

\medskip

\centerline{$\hat k(\alpha,\beta)=k^{(\alpha,\beta)}$}

\medskip

\noindent and 

\medskip

\centerline{$\hat m(\alpha,\beta,\gamma)=m^{(\alpha,\beta)}(\gamma)$}

\medskip

\noindent and for $t\leq\hat m(\alpha,\beta,\gamma)$ let 

\medskip

\centerline{$\hat\eta(\alpha,\beta,\gamma,t)=\eta^{(\alpha,\beta)}(\gamma,t)$.}

\medskip

Claim 2. Suppose $n\in\omega$. Then we have the following:

(1) $\hat\jmath(0,0,0)>y(n+\hat m(0,0,0)+1)$, and

(2) For every $\alpha\leq k^*$ and $\beta\leq m^*(\alpha)$ and $\gamma<\hat k(\alpha,\beta)$ we have
$\hat\jmath(\alpha,\beta,\gamma+1)>y(\hat\jmath(\alpha,\beta,\gamma)+\hat m(\alpha,\beta,\gamma+1)+1)$, and

(3) For every $\alpha\leq k^*$ and every $\beta< m^*(\alpha)$  we have
$\hat\jmath(\alpha,\beta+1,0)>y(\hat\jmath(\alpha,\beta,\allowbreak\hat k(\alpha,\beta))+\hat m(\alpha,\beta+1,0)+1)$, and

(4) For every $\alpha< k^*$ we have
$\hat\jmath(\alpha+1,0,0)>y(\hat\jmath(\alpha,m^*(\alpha),\hat k(\alpha,m^*(\alpha)))+\hat m(\alpha+1,0,0)+1)$.

Proof: Clause (1) holds because $j^{(0,0)}(0)>y(n+m^{(0,0)}(0)+1)$.

Clause (2) holds because $j^{(,\alpha,\beta)}(\gamma+1)>y(j^{(\alpha,\beta)}(\gamma)+m^{(\alpha,\beta)}(\gamma+1)+1)$.

We verify clause (3) as follows.

Case A: $\alpha=0$.

Notice that $j^{0,\beta)}(k^{(0,\beta)})<h(z(n+\beta+1))$ becuase the system $(j^{(0,\beta)},\allowbreak
m^{(0,\beta)},\allowbreak k^{(0,\beta)})$ is $h$-tight above $z(n+\beta+1)$. Notice also that
$j^{(0,\beta+1)}(0)>y(z(n+\beta+2)+m^{(0,\beta+1)}(0)+1)$ because the system $(j^{(0,\beta+1)},\allowbreak
m^{(0,\beta+1)},\allowbreak k^{(0,\beta+1)})$ is above $z(n+\beta+2)$.
 Hence we have
$\hat\jmath(0,\beta+1,0)=j^{(0,\beta+1)}(0)>y(z(n+\beta+2)+m^{(0,\beta+1)}(0)+1)\geq
y(h(z(n+\beta+1))+m^{(0,\beta+1)}(0)+1)\geq
y(j^{(0,\beta)}(k^{(0,\beta)})+m^{(0,\beta+1)}(0)+1)=y(\hat\jmath(0,\beta,\hat k(0,\beta))+\hat m(0,\beta+1,0)+1)$.

Case B: $\alpha>0$.

Notice that $j^{\alpha,\beta)}(k^{(\alpha,\beta)})<h(z(j^*(\alpha-1)+\beta+1))$ becuase the system $(j^{(\alpha,\beta)},\allowbreak
m^{(\alpha,\beta)},\allowbreak k^{(\alpha,\beta)})$ is $h$-tight above $z(j^*(\alpha-1)+\beta+1)$. Notice also that
$j^{(\alpha,\beta+1)}(0)>y(z(j^*(\alpha-1)+\beta+2)+m^{(\alpha,\beta+1)}(0)+1)$ because the system $(j^{(\alpha,\beta+1)},\allowbreak
m^{(\alpha,\beta+1)},\allowbreak k^{(\alpha,\beta+1)})$ is above $z(j^*(\alpha-1)+\beta+2)$.
 Hence we have
$\hat\jmath(\alpha,\beta+1,0)=j^{(\alpha,\beta+1)}(0)>y(z(j^*(\alpha-1)+\beta+2)+m^{(\alpha,\beta+1)}(0)+1)\geq
y(h(z(j^*(\alpha-1)+\beta+1))+m^{(\alpha,\beta+1)}(0)+1)\geq
y(j^{(\alpha,\beta)}(k^{(\alpha,\beta)})+m^{(\alpha,\beta+1)}(0)+1)=y(\hat\jmath(\alpha,\beta,\hat k(\alpha,\beta))+\hat m(\alpha,\beta+1,0)+1)$.

To see that clause (4) holds, we have $\hat\jmath(\alpha+1,0,0)=j^{(\alpha+1,0)}(0)\geq
y(z(j^*(\alpha)+1)+m^{(\alpha+1,0)}(0)+1)\geq y(h(z(j^*(\alpha)))+m^{(\alpha+1,0)}(0)+1)
\geq  y(h(j^*(\alpha))+m^{(\alpha+1,0)}(0)+1)
\geq y(h(z(j^*(\alpha-1))+m^*(\alpha)+1))+m^{(\alpha+1,0)}(0)+1)\geq
y(j^{(\alpha,m^*(\alpha))}(k^{(\alpha,m^*(\alpha))})+m^{(\alpha+1,0)}(0)+1)$.

The first inequality is because the system $(j^{(\alpha+1,0)},m^{(\alpha+1,0)},k^{(\alpha+1,0)})$ is above
$z(j^*(\alpha)+1)$ whence by clause (2) of Definition 7.1 we have the first inequality. 
The second inequality is by definition of the function $z$.  The third inequality is by the fact that $z$ is an increasing function. 
The fourth inequality is because $(j^*,m^*,k^*)$ satisfies clause (2) of Definition 7.1. The fifth inequallity
is because the system $(j^{(\alpha,m^*(\alpha))},m^{(\alpha,m^*(\alpha))},k^{(\alpha,m^*(\alpha))})$ is $h$-tight 
above $z(j^*(\alpha-1)+m^*(\alpha)+1)$.

The Claim is established.

Claim 3.  Suppose $n\in\omega$ and $\nu\in H$. Then there are $\alpha\leq k^*$ and $\beta\leq m^*(\alpha)$ and
$\gamma\leq\hat k(\alpha,\beta)$ and $\delta\leq \hat m(\alpha,\beta,\gamma)$ such that
$\nu$ is comparable with $\hat\eta(\alpha,\beta,\gamma,\delta)$.

Proof. Pick $t\in{\cal C}$ such that $\nu\in t$. Take  $\alpha$ and $\beta$ such that
$t$  is comparable with $\eta^*(\alpha,\beta)$. 

Case 1: $\nu\in\eta^*(\alpha,\beta)$.

Take $\gamma\leq k^{(\alpha,\beta)}$ and
$\delta\leq m^{(\alpha,\beta)}(\gamma)$ such that $\nu$ is comparable with $\eta^{(\alpha,\beta)}(\gamma,\delta)$.

Case 2: $\nu\notin\eta^*(\alpha,\beta)$.

If $\alpha = 0$ then let $\zeta= z(n+\beta+1)$ and if $\alpha>0$ then let
$\zeta= z(j^*(\alpha-1)+\beta+1)$. Let $\nu'=\nu\restr h(\zeta)$. Choose $\gamma\leq k^{(\alpha,\beta)}$ and
$\delta\leq m^{(\alpha,\beta)}(\gamma)$ such that $\nu'$ is comparable with
$\eta^{(\alpha,\beta)}(\gamma,\delta)$. Becuase the system $(j^{\alpha,\beta}),m^{\alpha,\beta)},k^{(\alpha,\beta)}$ is
$h$-tight above $\zeta$ we have $\eta^{(\alpha,\beta)}(\gamma,\delta)\leq \nu'$. Therefore
$\eta^{(\alpha,\beta)}(\gamma,\delta)\leq \nu$.

The Claim is established.

For each $n\in\omega$ define $\zeta(\alpha,\beta,\gamma)$ by the following recursive formulas:
$$\zeta(0,0,0)=0.$$
For $\alpha\leq k^*$ and $\beta\leq m^*(\alpha)$ and $\gamma<\hat k(\alpha,\beta)$ we have
$$\zeta(\alpha,\beta,\gamma+1)=\zeta(\alpha,\beta,\gamma)+1.$$
For $\alpha\leq k^*$ and $\beta< m^*(\alpha)$ we have
$$\zeta(\alpha,\beta+1,0)=\zeta(\alpha,\beta,\hat k(\alpha,\beta))+1.$$
For $\alpha< k^*$ we have
$$\zeta(\alpha+1,0,0)=\zeta(\alpha,m^*(\alpha),\hat k(\alpha,m^*(\alpha)))+1.$$

Define $\tilde\jmath(\zeta(\alpha,\beta,\gamma))=\hat\jmath(\alpha,\beta,\gamma)$, and
 $\tilde m(\zeta(\alpha,\beta,\gamma))=\hat m(\alpha,\beta,\gamma)$, and
 $\tilde k(\zeta(\alpha,\beta))=\hat k(\alpha,\beta)$, and
 $\tilde\eta(\zeta(\alpha,\beta,\gamma,\delta))=\hat\eta(\alpha,\beta,\gamma,\delta)$, and

Claim 4. $(\tilde \jmath, \tilde k,\tilde m)$ is a $y$-bound system above $n$ and $H$ is a 
$(\tilde\jmath,\tilde k,\tilde m,\tilde\eta)$-squeezed tree.

Proof. By Claims 2 and 3.

The Lemma is established.

\proclaim Lemma 7.13.  Suppose $x\in{}^\omega(\omega-\{0\})$ is strictly increasing and suppose
that for each $n\in\omega$ we have that  $T_n$ is an $x$-squeezed tree.  Then there are
 $T^*$ and $\langle\gamma_t\,\colon\allowbreak
t\in\omega\rangle$ an increasing sequence of integers such that
$T^*$ is an $x$-squeezed tree and $\gamma_0=0$ and
$(\forall t>0)\allowbreak(t<\gamma_t)$ and for every 
$f\in{}^{<\omega}\omega$ we have

\centerline{$(\forall t>0)(\exists s< t)(f\restr\gamma_t\in T_{\gamma_s})$ iff $f\in T^*$.}

\medskip

Proof: For each $n\in\omega$ choose $h_n\in{}^\omega\omega$ such that $T_n$ is $h_n$-tight.

We build as follows. Let $\gamma_0=0$.  Given $\gamma_t$, define $g_t(0)=\gamma_t$. For $0\leq s\leq t$ let
$g_t(s+1)=h_{\gamma_s}(g_t(s))$. Let $\gamma_{t+1}=g_t(t+1)$.

Let $T^*=\{\eta\in{}^{<\omega}\omega\,\colon\allowbreak (\forall t>0)(\exists s< t)
(\eta\restr \gamma_t\in T_{\gamma_s})\}$.

Now fix $n\in\omega$.  We build an $x$-bound system $(j,m,k)$ above $n$ and we build $\eta$ so that
$(j,m,k)$ and $\eta$ witness the fact that $T^*$ is $x$-squeezed.

For every $t\in\omega$ and $s\leq t$ choose an $h_{\gamma_s}$-tight $x$-bound system $(j^s_t,m^s_t,k^s_t)$ above
$g_t(s)$ along with $\eta^s_t$ such that
$T_{\gamma_s}$ is $(j^s_t,m^s_t,k^s_t,\eta^s_t)$-squeezed.

We define $\zeta$ such that for $\alpha\geq n$ and $\beta\leq \alpha$ and $\gamma\leq k^\alpha_\beta$ we have

\begin{itemize}

\item $\zeta(n,0,0)=0$, and

\item if $\gamma < k^\alpha_\beta$ then $\zeta(\alpha,\beta,\gamma+1)=\zeta(\alpha,\beta,\gamma)+1$, and

\item if $\beta< \alpha$ then $\zeta(\alpha,\beta+1,0)=\zeta(\alpha,\beta, k^\alpha_\beta)+1$, and

\item  if $\alpha\geq n$ then $\zeta(\alpha+1,0,0)=\zeta(\alpha,\alpha,k^\alpha_\alpha)+1$.

\end{itemize}

We define $(j,m,k)$ such that for every $\alpha\geq n$ and $\beta\leq \alpha$ and $\gamma\leq k^\alpha_\beta$ we have

\begin{itemize}

\item $j(\zeta(\alpha,\beta,\gamma))=j^\alpha_\beta(\gamma)$, and

\item $m(\zeta(\alpha,\beta,\gamma))=m^\alpha_\beta(\gamma)$, and

\item $k=k^\alpha_\beta$.

\end{itemize}

Claim 1. $(j,m,k)$ is an $x$-bound system above $n$.

Proof: Clause (1) of Definition 7.1 is immediate.

Clause (2) of Definition 7.1 holds because $j(0)=j^0_n(0)>x(g^0_n(0)+m^0_n(0)+1)\geq x(n+m(0)+1)$.  The first inequality holds because 
the system $(j^0_n,m^0_n,k^0_n)$ is above $g^0_n(0)$ and it satisfies clause (2) of Definition 7.1.

We have $j(\zeta(\alpha,\beta,\gamma+1))=j^\beta_\alpha(\gamma+1)>x(j^\beta_\alpha(\gamma)+m^\beta_\alpha(\gamma+1)+1)=
x(j(\zeta(\alpha,\beta,\gamma))+m(\zeta(\alpha,\beta,\gamma+1))+1)$.

We have $j(\zeta(\alpha,\beta+1,0))=j^{\beta+1}_\alpha(0)>x(g_\alpha(\beta+1)+m^{\beta+1}_\alpha(0)+1)\geq
x(h_{\gamma_\beta}(g_\alpha(\beta))+m^{\beta+1}_\alpha(0)+1)\geq
x(j^\beta_\alpha(k^\beta_\alpha)+m^{\beta+1}_\alpha(0)+1)=x(j(\zeta(\alpha,\beta,k^\beta_\alpha))+m(\zeta(\alpha,\beta+1,0))+1)$.

The first inequality is clause (2) of Definition 7.1 applied to the system $(j^{\beta+1}_\alpha,\allowbreak
m^{\beta+1}_\alpha,k^{\beta+1}_\alpha)$.  The second inequality is by the definition
 of $g_\alpha$. The third inequality is because the system
$(j^\beta_\alpha,m^\beta_\alpha,k^\beta_\alpha)$ is $h_{\gamma_\beta}$-tight above $g_\alpha(\beta)$.

The Claim is established.

We define $\eta$ such that for every $\alpha\geq n$ and $\beta\leq \alpha$ and $\gamma\leq k^\alpha_\beta$ and $\delta\leq m^\beta_\alpha(\gamma)$
we have $\eta(\zeta(\alpha,\beta,\gamma),\delta)=\eta^\beta_\alpha(\gamma,\delta)$.

Claim 3: $T^*$  is a $(j,k,m,\eta)$-squeezed tree. 

Proof: It is straightforward to verify that $T^*$ is a tree and that clause (1) and clause (2) of Definition 7.2 hold.

To verify clause (3), suppose we have $\nu\in T^*$. We show that $\nu$ is comparable to some $\eta(l,i)$
with $(l,i)\in{\rm dom}(\eta)$. Choose $\nu'\in T^*$ such that $\nu\leq \nu'$ and ${\rm lh}(\nu')\geq
\gamma_{n+1}$. It suffices to show that $\nu'$ is comparable with some $\eta(l,i)$ with $(l,i)\in {\rm dom}(\eta)$.
Because $\nu'\in T^*$ we may choose $s\leq n$ such that $\nu'\restr\gamma_{n+1}\in T_{\gamma_s}$.
We may select $(l,i)\in{\rm dom}(\eta^s_n)$ such that $\nu'\restr\gamma_{n+1}$ is comparable with $\eta^s_n(l,i)$.
We have $\eta^n_s(l,i)=\eta(\zeta(n,s,l),i)$, so ${\rm lh}(\eta^s_n(l,i))=j(\zeta(n,s,l))=j^s_n(l)\leq
h_{\gamma_s}(g_n(s))=g_n(s+1)\leq\gamma_{n+1}$. Therefore $\eta(\zeta(n,s,l),i)\leq\nu'\restr\gamma_{n+1}$ and therefore
$\eta(\zeta(n,s,l),i)$ is comparable with $\nu'$.

The Claim and the Lemma are established.

\proclaim Lemma 7.14.  Suppose $x\in{}^\omega(\omega-\{0\})$ is strictly increasing, and suppose
 $P$ is a forcing notion such that\/ {\rm $V[G_P]\models``$for all countable $X\subseteq V$ there is
a countable $Y\in V$ such that $X\subseteq Y$ and $\langle T_n\,\colon\allowbreak
n\in\omega\rangle$ is a sequence of  $x$-squeezed trees and $(\forall n\in\omega)\allowbreak(T_n\in V)$.''}
Then\/ {\rm $V[G_P]\models``$there is a strictly increasing sequence of integers $\langle m_i\,\colon\allowbreak
i\in\omega\rangle$ and an $x$-squeezed tree $T^*\in V$ such that $m_0=0$ and $(\forall i>0)\allowbreak(m_i>i)$ and
for every $\eta\in{}^{<\omega}\omega$, if
$(\forall i>0)\allowbreak(\exists j<i)\allowbreak(\eta\restr m_{i+1}\in T_{m_j})$ then $\eta\in T^*$.''}

Proof: Work in $V[G_P]$. Let $b\in V$ be a countable set such that
$\{T_n\,\colon\allowbreak n\in\omega\}\subseteq b\in V$ and $(\forall x\in b)\allowbreak(x$ is an $x$-squeezed tree).
Let $\langle S_n\,\colon\allowbreak n\in\omega\rangle\in V$ enumerate $b$ with infinitely many repetitions such that
$S_0=T_0$.
Build $\langle S'_n\,\colon\allowbreak n\in\omega\rangle$ by setting $S'_0=S_0$ and for every $n>0$ set
$S'_n=S_n\cup S'_{n-1}$. Build $h$ mapping $\omega$ into $\omega$ inductively by setting $h(0)=0$ and for every
$n>0$ set $h(n)$ equal to the least integer $m$ such that $m>h(n-1)$ and $T_n=S_m$.

Using Lemma 7.13, take $ T^*\in V$  an $x$-squeezed tree and $\langle k_i\,\colon\allowbreak i\in\omega\rangle\in V$ such that
for every $\eta\in{}^{<\omega}\omega$ we have $\eta\in T^*$ iff
$(\forall n>0)\allowbreak(\exists i<n)\allowbreak(\eta\restr k_n\in S'_{k_i})$.

Build $\langle n'_i\,\colon\allowbreak i\in\omega\rangle$ an increasing sequence of integers such that
$n'_0=0$ and $n'_1>k_1$ and for every $i\in\omega$ we have $h(n'_i)<n'_{i+1}$ and 

(*)  $(\exists t\in\omega)\allowbreak
(n'_i<k_t<n'_{i+1})$.

For every $i\in\omega$ let $m_i=h(n'_{3i+3})$.

Fix $\eta\in{}^{<\omega}\omega$ such that $(\forall i>0)(\exists j<i)(\eta\restr m_{i+1}\in T_{m_j})$. To establish the Lemma, it suffices to show $\eta\in T^*$.
By choice of $T^*$, it suffices to show $(\forall n>0)\allowbreak(\exists i<n)\allowbreak(\eta\restr k_{n}\in S'_{k_i})$.

Claim 1. $(\forall i>0)(\exists j<i)(\eta\restr n'_{i+1}\in S'_{n'_j})$.

Proof: The proof breaks into two cases.

Case 1: $i<6$.

We have $n'_{i+1}\leq n'_6\leq h(n'_6)= m_1$, and $\eta\restr m_1\in T_0$, so
$\eta\restr n'_{i+1}\in S_0'$.

Case 2: $i\geq 6$.

Fix $i^*>0$ such that $3i^*+3\leq i \leq 3i^*+5$.

We may fix $j^*<i^*$ such that $\eta\restr m_{i^*+1}\in T_{m_{j^*}}$.  

Now, we have 

(*)  $i+1\leq 3i^*+6$ so

(**) $ n'_{i+1}\leq h(n'_{3i^*+6})= m_{i^*+1}$. 

 We also have 

(***) $\eta\restr m_{i^*+1}\in T_{m_{j^*}}\subseteq S'_{h(m_{j^*})}$.

By (**) and (***) we have

(****)  $\eta\restr n'_{i+1}\in S'_{h(m_{j^*})}$.

Note that

(*****) $h(m_{j^*)}=h(h(n'_{3j^*+3}))\leq n'_{3j^*+5}\leq n'_{3i^*+2}\leq n'_{i-1}$.

By (****) and (*****) we have $\eta\restr n'_{i+1}\in S'_{h(m_{j^*})}\subseteq S'_{n'_{i-1}}$.

The Claim is established.

To complete the proof of the Lemma, suppose $i>0$. We must show that there is $t<i$ such that
$\eta\restr k_{i}\in S'_{k_t}$.

Case 1: $k_{i-1}<n'_0$.

By (*) we have $n'_1\geq k_{i}$. By Claim 1 we have $\eta\restr n'_1\in S_0$.  Hence
$\eta\restr k_{i}\in S_0$.

Case 2:  $n'_0\leq k_{i-1}$.

By (*) we know that there is at most one element of $\{n'_j\,\colon\allowbreak j\in\omega\}$ strictly between
$k_{i-1}$ and $k_{i}$. Hence we may
fix $j>0$ such that $n'_{j-1}\leq k_{i-1}<k_{i}\leq n'_{j+1}$.  If $\eta\restr n'_{j+1}\in S_0$ then
$\eta\restr k_{i}\in S_0$ and we are done, so assume otherwise.  By Claim 1 we may fix $m<j$ such that $\eta\restr n'_{j+1}\in S'_{n'_m}$.
We have $\eta\restr k_{i}\in S'_{n'_m}\subseteq S'_{n'_{j-1}}\subseteq S'_{k_{i-1}}$ and again we are done.

The Lemma is established.

\proclaim Lemma 7.15.  Suppose  $P$ is a forcing notion and
$y\in{}^\omega(\omega-\{0\})$ is strictly increasing. Suppose  
 $\langle T_n\,\colon\allowbreak n\in\omega\rangle$ is a sequence of $y$-squeezed trees.  Then there is a $y$-squeezed tree
$T^*$ such that in $V[G_P]$ we have that for every
$n\in\omega$ and every $j\in\omega$ and every $g\in[T_j]$
there is $k\in\omega$ such that for every $\eta\in T_j$ extending $g\restr k$, if 
 $\eta\restr k\in  T^*$ then
$\eta\in  T^*$.

Proof:    Build a sequence of $y$-squeezed trees  $\langle T'_j\,\colon j\in\omega\rangle$ such that $T'_0=T_0$ and for every
$j\in\omega$ we have $T'_{j+1}=T'_j\cup T_{j+1}$. 
By Lemma 7.13 we may find an increasing sequence of integers 
$\langle k_n\,\colon\allowbreak n\in\omega\rangle$ and a $y$-squeezed tree $T^*$ such that $k_0=0$ and
$(\forall n>0)\allowbreak(k_n>n)$ and
for every $\eta\in{}^{<\omega}\omega$ we have

\medskip

{\centerline{$(\forall n>0)(\exists i< n)(
\eta\restr k_n\in T'_{k_i})$ iff $\eta\in T^*$.}}

\medskip

Fix a forcing notion $P$ and work in $V[G_P]$. Fix  $j\in\omega$ and $g\in [T_j]$.
Let $k={\rm max}\{k_{j'}\colon\allowbreak j'\leq j\}$.  Fix $\eta\in T_j$ extending $g\restr k$ and assume
$\eta\restr k\in T^*$.  It suffices to show that $\eta\in T^*$. If $j=0$ then $\eta\in T_0=T'_0\subseteq T^*$.
Therefore, we assume that $j>0$.  It suffices to show that

\medskip

\centerline{ $(\forall i>0)\allowbreak(\exists i'< i)
\allowbreak(\eta\restr k_i\in T'_{k_{i'}})$. }

\medskip

 Towards this end, fix $i>0$.

Case 1: $i\leq j$.

Because $k_i\leq k$ we have $\eta\restr k_i\in T^*$. Therefore we may take $i'< i$ such that
$\eta\restr k_{i}\in T'_{k_{i'}}$.  

Case 2: $0<j<i$.

Because $\eta\in T_j$ we  have $\eta\restr i\in T_j\subseteq T'_{j+1}\subseteq T'_{k_j}$.

The Lemma is established.

\proclaim Lemma 7.16.  Suppose   $y\in{}^\omega(\omega-\{0\})$ is strictly increasing.
Suppose that for every $n\in\omega$ we have
that
$T_n$ is a $y$-squeezed tree. 
Suppose  $\zeta\in{}^\omega\omega$.
Then there is a $y$-squeezed tree $T^*\supseteq T$ and a sequence of integers
$\langle m_i\,\colon\allowbreak i\in\omega\rangle$
such that $\zeta\in[T^*]$ and
for every  $i\in\omega$ and every $j>m_i$ and every $\nu\in T_{m_i}$ extending $\zeta\restr j$
we have  $\nu\in T^*$.

Proof.  Define $\langle T'_k\,\colon\allowbreak k\in\omega\rangle $ by setting $T'_0=
 T_0\cup\{\zeta\restr n\,\colon\allowbreak n\in\omega\}$ and for every $k\in\omega$
set $T'_{k+1}= T'_k\cup T_{k+1}$.

By Lemma 7.13 we may choose $T^*$ a $y$-squeezed tree and
$\langle m_i\,\colon\allowbreak i\in\omega\rangle$ an increasing sequence of integers
such that

\medskip

{\centerline{$(\forall g\in{}^\omega\omega)\allowbreak((\forall n>0)\allowbreak
(\exists i< n)\allowbreak(g\restr m_n\in T'_{m_i})$ iff $g\in[T^*])$.}}

\medskip

Now suppose that $\eta\in T$ and $i\in\omega$ and
${\rm length}(\eta)\geq m_i$ and
 $\nu$ extends $\eta$
and $\nu\in T_{m_i}$.  We show $\nu\in T^*$.

Choose $h\in[T_{m_i}]$ extending $\nu$.  It suffices to show that $h\in[T^*]$.
Therefore it suffices to show that $(\forall k>0)\allowbreak(\exists
j< k)\allowbreak(h\restr m_k\in T'_{m_j})$.

Fix $k\in\omega$. 
If $i< k$ then because $h\in[T_{m_i}]$ we have that
$h\restr m_k\in T_{m_i}\subseteq T'_{m_i}$ and we are done.  If instead $k\leq i$ then
$h\restr m_k=\eta\restr m_k\in T'_{0}$ and again we are done.

The Lemma is established.

\proclaim Theorem 7.17.  Suppose $\langle P_\eta\,\colon\allowbreak\eta\leq\kappa\rangle$ is a countable support iteration based on $\langle Q_\eta\,\colon\allowbreak\eta<\kappa\rangle$ and suppose\/ {\rm $(\forall\eta<\kappa)\allowbreak({\bf 1}\forces_{P_\eta}``Q_\eta$ 
is proper and has the $P$-point property'').} 
Suppose $\lambda$ is a sufficiently large regular cardinal 
and $\alpha<\kappa$ 
and  $y$ is a $P_\alpha$-name and $T$ is a $P_\kappa$-name 
and\/ {\rm ${\bf 1}\forces_{P_\alpha}``y\in
{}^\omega(\omega-\{0\})$ is strictly increasing''} and\/ {\rm
 ${\bf 1}\forces_{P_\kappa}``T$ is a $y$-squeezed tree.''}
Suppose $N$ is a countable elementary submodel of $H_\lambda$ and $\{P_\kappa, \alpha,  y, T
\}\in N$. 
Suppose
 $p\in P_\alpha$ and $p$ is $N$-generic.
Then\/ {\rm $p\forces``(\forall q\in P_{\alpha,\kappa})\allowbreak(\exists q'\leq q)\allowbreak(\exists H)\allowbreak(H$ is
a $y$-squeezed tree and $q'\forces`T\subseteq H$').''}

Proof: The proof proceeds by induction on $\kappa$.  We assume that 
  $\lambda$, 
$N$, $\alpha$, $p$,  $y$, and $T$ are as in the hypothesis of the Theorem.
Fix a $P_\alpha$-name $q$ in $N$ such that ${\bf 1}\forces``q\in P_{\alpha,\kappa}$.

Case 1. $\kappa=\beta+1$.

Using Lemma 7.12, choose $\tilde q$ and $H'$ such that 

\medskip

\centerline{
${\bf 1}\forces_{P_\beta}``\tilde q\leq q(\beta)$ and $H'$ is a $y$-squeezed tree and
$\tilde q\forces``T\subseteq H'$.'\thinspace''}

\medskip

We may assume that the names $\tilde q$ and $H'$ are elements of $N$.
Use the induction hypothesis to take a $P_\alpha$-names $q^*$ and
 $H$  such that

\medskip

\centerline{
$p\forces``q^*\leq q\restr\beta$ and $H$ is a $y$-squeezed tree and
$q^*\forces`H'\subseteq H$.'\thinspace''}

\medskip

We have that $p\forces``(q^*,\tilde q)\forces`T\subseteq H$.'\thinspace'''
Case 1 is established.

Case 2.  ${\rm cf}(\kappa)>\omega$.

Because no $\omega$-sequences of ordinals can be added at limit stages of uncountable cofinality, 
we may take $\beta$ and $T'$ and $q'$ to be $P_\alpha$-names in $N$ such that

\medskip

\centerline{
${\bf 1}\forces``\alpha\leq\beta<\kappa$ and $T'$ is a
$P_{\alpha,\beta}$-name and $q'\leq q$ and}

\centerline{
${\bf 1}\forces_{P_{\alpha,\beta}}`T'$ is a $y$-squeezed tree' and
$q'\forces_{P_{\alpha,\kappa}}``T'=T$.''\thinspace'\thinspace''}

\medskip

For every $\beta_0\in\kappa\cap N$ such that $\alpha\leq\beta_0$ 
let $\tilde q(\beta_0)$ and $H(\beta_0)$  be  $P_\alpha$-names in $N$ such that

\medskip

\centerline{
${\bf 1}\forces``$if $\beta=\beta_0$ and there is some $\tilde q\leq q'\restr\beta$ and some 
$H^*$}

\centerline{such that $H^*$ is a $y$-squeezed tree and
$\tilde q\forces`T'\subseteq H^*$,'}

\centerline{then $q^*(\beta_0)$ and $H(\beta_0)$  are witnesses thereto.''}

\medskip

Let $ q^*$ and $H$  and $s$ be $P_\alpha$-names such that for every $\beta_0
\in \kappa\cap N$, if $\alpha\leq\beta_0$, then

\medskip

\centerline{
${\bf 1}\forces``$if $\beta=\beta_0$ then $ q^*=q^*(\beta_0)$ and $H=H(\beta_0)$
 and 
$s\in P_{\alpha,\kappa}$ and}

\centerline{$s\restr\beta=q^*$ and
$s\restr[\beta,\kappa)=q'\restr[\beta,\kappa)$.''}

\medskip

Claim 1: $p\forces``s\leq q$ and $s
\forces`T\subseteq H$.'\thinspace''

Proof:  Suppose $p'\leq p$. Choose $p^*\leq p'$ and $\beta_0<\kappa$ such that
$p^*\forces``\beta_0=\beta$.''  Because the name $\beta$ is in $N$ and $p^*$ is $N$-generic,
we have that $\beta_0\in N$.  Notice by the induction hypothesis
we have 

\medskip

\centerline{$p\forces``$there is some $q^\#\leq q'\restr\beta_0$ and 
some $y$-squeezed tree
$H^\#$ }

\centerline{
such that $q^\#\forces`T'\subseteq H^\#$.'\thinspace''}

\medskip

Hence 

\medskip

\centerline{$p^*\forces``q^*=q^*(\beta_0)\leq q'\restr\beta$ and $H=H(\beta_0)$ and
$H$ is a $y$-squeezed tree and}

\centerline{ $q^*\forces`T'\subseteq H$ and
$q'\restr[\beta,\kappa)\forces``T'=T$.''\thinspace'\thinspace'' }

\medskip

 Therefore
$p^*\forces``s\forces`T\subseteq H$.'\thinspace''

Claim 1 is established.
This completes Case 2.

Case 3. ${\rm cf}(\kappa)=\omega$.

Let $\langle\alpha_n\,\colon\allowbreak n\in\omega\rangle$ be an increasing 
sequence from $\kappa\cap N$ cofinal in $\kappa$ such that 
$\alpha_0=\alpha$.

By Theorem 3.3 we have that ${\bf 1}\forces_{P_\alpha}``P_{\alpha,\kappa}$ is $\omega^\omega$-bounding,'' so we may fix $P_\alpha$-names
$\tilde q$ and $h$ 
 such that ${\bf 1}\forces_{P_\alpha}``h\in{}^\omega\omega$ and $\tilde q\leq q$ and $\tilde q\forces_{P_{\alpha,\kappa}}`
T$ is $h$-tight.'\thinspace'' We may assume the names $h$ and $\tilde q$ are in $N$.

Working in $V[G_{P_{\alpha}}]$, let $z$ and
 $\langle{\cal T}_\beta\,\colon\beta\in\omega\rangle$ and ${\cal T}$ be as in the proof of Lemma 7.12.
In $V[G_{P_\kappa}]$ let $\zeta\in[{\cal T}]$ be defined by $(\forall n\in\omega)\allowbreak(\zeta(n)=T\cap{}^{<h(n)}\omega)$.
In $V[G_{P_\alpha}]$ fix an isomorphism from ${}^{<\omega}\omega$ onto ${\cal T}$ and implicitly fix a $P_\alpha$-name for the isomorphism that is
an element of $N$.

Using Lemma 3.1, fix
 $\langle (p_n,\zeta_n)\,\colon\allowbreak n\in\omega\rangle\in N$ (that is, the sequence of names is an element of $N$ but not necessarily their values) such that 
${\bf 1}\forces``p_0\leq \tilde q$'' and 
 for every $n\in\omega$ we have that each of the following holds:

(0) $p_n$ is a $P_\alpha$-name for an element of $P_{\alpha,\kappa}$, and

(1) For every $k\leq n$ we have ${\bf 1}\forces_{P_{\alpha_n}}``p_0\restr[\alpha_n,\kappa)
\forces`\zeta\restr k=\zeta_n\restr k$,'\thinspace'' and

(2) $\zeta_n$ is a $P_{\alpha_n}$-name for an element of $[{\cal T}]$, and

(3) ${\bf 1}\forces_{P_{\alpha_n}}``p_0\restr[\alpha_n,\alpha_{n+1})\forces`\zeta_n\restr k
=\zeta_{n+1}\restr k$ for every $k\leq n+1$,'\thinspace'' and

(4) ${\bf 1}\forces_{P_\alpha}``p_{n+1}\leq p_n$,'' and

(5) whenever $k\leq m<\omega$ we have ${\bf 1}\forces_{P_{\alpha_n}}``p_m\restr[
\alpha_n,\alpha_{n+1})\forces`\zeta_n\restr k=\zeta_{n+1}\restr k$.'\thinspace''

Claim 2.  Suppose $\alpha\leq\beta\leq\gamma<\kappa$ and 
suppose $T'$ is a $P_\gamma$-name for a $z$-squeezed tree. 
Then

\medskip

\centerline{ ${\bf 1}\forces_{P_\beta}``V[G_{P_\alpha}]\models`(\forall q\in P_{\beta,\gamma})(\exists q'\leq q)(\exists H)$}

\centerline{($H$ is a $z$-squeezed tree and $q'\forces``T'\subseteq H$'').'\thinspace''}

\medskip

Proof: Given $r_1\in P_\alpha$ and a $P_\alpha$-name $r_2$ for an element of $P_{\alpha,\beta}$ and a 
$P_\beta$-name $q$ for an element of $P_{\beta,\gamma}\cap V[G_{P_\alpha}]$, choose $\lambda'$ a sufficiently large
regular cardinal and $N'$ a countable elementary substructure of $H_{\lambda'}$ containing $\{r_1,r_2,q,P_\kappa,\alpha,\beta,\gamma,z, T'\}$.
Choose $r_1'\leq r_1$ such that $r_1'$ is $N'$-generic.  By the overall induction hypothesis (i.e., because $\gamma<\kappa$)
we may choose $s$ such that

\medskip

\centerline{$r'_1\forces`` s\leq(r_2,q)$ and $(\exists H)(H$ is a $z$-squeezed tree and $s\forces`T'\subseteq H$').''}

\medskip

Consequently we may choose $H$ such that

\medskip

\centerline{$(r'_1,s\restr\beta)\forces``V[G_{P_\alpha}]\models`
s\restr[\beta,\gamma)\leq q$ and $H$ is a $z$-squeezed tree and}

\centerline{$s\restr[\beta,\gamma)\forces``T'\subseteq H$'').'\thinspace''}

\medskip

The Claim is established.

Let $\langle T'_{j}\,\colon\allowbreak j\in\omega\rangle$ list all $P_{\alpha}$-names $T'\in N$ such that 
 we have that $T'$ is a
$z$-squeezed tree.

Using Lemma 7.15, choose $T^*\supseteq T$ a $z$-squeezed tree such that, in $V[G_{P_\kappa}]$, we have that
for every $n\in\omega$ and every $j\in\omega$ and every $g\in[T'_j]$ there exists
$k\in\omega$ such that for every $\eta\in T'_j$ extending $g\restr k$, if
$\eta\restr k\in T^*$ then $\eta\in T^*$.

In the prceding paragraph, we worked in $V[G_{P_\kappa}]$ so that the brackets about $T'_j$ would be interpreted in $V[G_{P_\kappa}]$; i.e.,
$g$ net not be in $V[G_{P_\alpha}]$.

Claim 3.  We may be build $\langle r_n\,\colon\allowbreak n\in\omega\rangle$ such that
$r_0=p$  and for every $n\in\omega$ we have that the following hold:

(1) $r_n\in P_{\alpha_n}$ is $N$-generic, and

(2) $r_{n+1}\restr\alpha_n=r_n$, and

(3) $r_n\forces``\zeta_n\in[T^*]$,'' and

(4) $p\forces``r_n\restr[\alpha,\alpha_n)\leq p_0\restr\alpha_n$.''

Proof: By induction on $n$.  For $n=0$ we have nothing to prove. Suppose we have $r_n$.

Usimg Claim 2, let $F_0$ and $F_2$  be $P_{\alpha_n}$-names  such that

\medskip

(*) {
 ${\bf 1}\forces``F_0$ and
$F_2$ are functions, both of which are in $V[G_{P_\alpha}]$,}

\centerline{and each of whose domains
 is equal to $P_{\alpha_n,\alpha_{n+1}}$, such that}

\centerline{$(\forall q'\in P_{\alpha_n,\alpha_{n+1}}\cap V[G_{P_\alpha}]
 )\allowbreak (F_0(q')\subseteq{\cal T}$ is a $z$-squeezed tree}

\centerline{and
$F_2(q')\leq q'$}

\centerline{and $F_2(q')\forces`\zeta_{n+1}\in[ F_0(q')]$').''}

\medskip

 We may assume that the names $F_0$ and
$F_2$  are in $N$. Notice that $F_0$ and $F_2$ depend on $n$, although this dependence is suppressed in our notation.

Working in $V[G_{P_{\alpha_n}}]$, use Lemma 7.16 to choose 
$\tilde T_n\subseteq{\cal T}$ a $z$-squeezed tree and
$\langle k_i\,\colon\allowbreak i\in\omega\rangle$ an increasing sequence of
integers (this sequence depends on $n$ but this fact is suppressed in our notation)
such that 

\medskip

 {\bf if}  $F_0$ is a function in $V[G_{P_\alpha}]$
  whose domains
 is equal to $P_{\alpha_n,\alpha_{n+1}}$, such that
$(\forall q'\in P_{\alpha_n,\alpha_{n+1}}\cap V[G_{P_\alpha}]
 )\allowbreak (F_0(q')\subseteq{\cal T}$ is a $z$-squeezed tree)

\medskip

{\bf then} $\zeta_n\in[\tilde T_n]$ and $\tilde T_n\in V[G_{P_\alpha}]$
and for every $\eta$  and every $i\in\omega$ and every
$\nu\in F_0(p_{k_i}\restr[\alpha_n,\alpha_{n+1}))$,
 if 
$\eta$ is a proper initial segment of $\zeta_n$ and ${\rm length}(\eta)\geq  k_i$ and 
$\nu$ extends $\eta$, then $\nu\in\tilde T_n$. 

\medskip

We may assume the $P_{\alpha_n}$-name $\tilde T_n$ is in $N$.

  By (*)  we have that

\medskip

\centerline{
$r_n\forces``\zeta_n\in[\tilde T_n]$ and for every $\eta$ and every $i\in\omega$}

\centerline{and every
$\nu\in F_0(p_{k_i}\restr[\alpha_n,\alpha_{n+1}))$,}

\centerline{
 if $\eta$ is a proper initial segment of $\zeta_n$ and}

\centerline{
${\rm length}(\eta)\geq k_i$ and 
$\nu$ extends $\eta$ then $\nu\in\tilde T_n$.''}

\medskip

Because  $\tilde T_n$ is a $P_{\alpha_n}$-name in $N$ forced to be in $V[G_{P_\alpha}]$, we conclude that by the $N$-genericity of $r_n$ that

\medskip

{\centerline{$r_n\forces``\tilde T_n\in N[G_{P_\alpha}]$.''}}

\medskip

Therefore there is a $P_{\alpha_n}$-name $m$ such that

\medskip

{\centerline{$r_n\forces``\tilde T_n=T_m'$.''}}

\medskip

Because $T^*$ was chosen as in the conclusion of Lemma 7.15, we may choose
 $k$ to be a $P_{\alpha_n}$-name for an integer such that

\medskip

(**) $r_n\forces``(\forall\eta\in \tilde T_n)($if $\eta$ extends $\zeta_n\restr k$ and $\eta\restr k\in  T^*$ 
then $\eta\in  T^*)$.''

\medskip

Choose $j$ and $K$ to be $P_{\alpha_n}$-names for integers such that
$r_n\forces``K=k_j\geq k$.''

Subclaim 1. 
$r_n\forces``F_2(p_K\restr[\alpha_n,\alpha_{n+1}))\forces
`\zeta_{n+1}\in[\tilde T_n]$.'\thinspace''

Proof.  It suffices to show

\medskip 

$r_n\forces``F_2(p_K\restr[\alpha_n,\alpha_{n+1}))\forces
`(\forall \nu)\allowbreak($if $\nu$ is a proper initial segment of $\zeta_{n+1}$ and ${\rm lh}(\nu)\geq K$ then
$\nu\in\tilde T_n)$.'\thinspace''

\medskip

Fix $\nu$ and $\eta$ such that

\medskip

$r_n\forces``F_2(p_K\restr[\alpha_n,\alpha_{n+1}))\forces
` \nu\in T_{n+1}$ and ${\rm lh}(\nu)\geq K$ and $\eta=\nu\restr K$.'\thinspace''

\medskip

By the definition of $\langle p_i\,\colon\allowbreak i\in\omega\rangle$ we have

\medskip

(***) $r_n\forces``p_K\restr[\alpha_n,\alpha_{n+1})\forces`\eta$ is an initial segment of $\zeta_{n+1}$.'\thinspace''

\medskip

By (*) we have

\medskip

(****) $r_n\forces``F_2(p_K\restr[\alpha_n,\alpha_{n+1}))\forces`\nu$ is an initial segment of $\zeta_{n+1}$ and
$\zeta_{n+1}\in [F_0(p_K\restr[\alpha_n,\alpha_{n+1}))]$.'\thinspace''

\medskip

Combining (***), (****), the definition of $\tilde T_n$ and the fact that
$r_n\forces``K\in\{k_i\,\colon\allowbreak i\in\omega\}$,'' we have that

\medskip

{\centerline{$r_n\forces``F_2(p_K\restr[\alpha_n,\alpha_{n+1}))\forces`\nu\in \tilde T_n$.'\thinspace''}}

The Subclaim is established.

Subclaim 2. $r_n\forces``F_2(p_K\restr[\alpha_n,\alpha_{n+1}))\forces`\zeta_{n+1}\in[ T^*]$.'\thinspace''

Proof: By (**) we have

\medskip

($\dag$) $r_n\forces``(\forall\eta\in \tilde T_n)\allowbreak(\eta\restr K\in  T^*$ implies
$\eta\in T^*)$.''

\medskip

Work in $V[G_{P_{\alpha_n}}]$ with $r_n\in G_{P_{\alpha_n}}$.
Fix $\eta\in\tilde T_n$ and suppose $F_2(p_K\restr[\alpha_n,\alpha_{n+1}))\forces``
\eta$ is an initial segment of $\zeta_{n+1}$ and ${\rm lh}(\eta)\geq K$.''
To establish the Subclaim it suffices to show

($\#)$ $F_2(p_K\restr[\alpha_n,\alpha_{n+1}))\forces``\eta\in T^*$.''

By the definition of $\langle p_i\,\colon\allowbreak i\in\omega\rangle$ we have

\medskip

{\centerline{$p_K\restr[\alpha_n,\alpha_{n+1})\forces``\eta\restr K=\zeta_n\restr K$.''}}

\medskip

Hence by the fact that Claim 3 holds for the integer $n$ we have

\medskip

($\dag\dag$) $p_K\restr[\alpha_n,\alpha_{n+1}))\forces``\eta\restr K\in  T^*$.''

\medskip

By Subclaim 1, ($\dag$), ($\dag\dag$), and the fact that
$F_2(p_K\restr[\alpha_n,\alpha_{n+1}))\leq
p_K\restr[\alpha_n,\alpha_{n+1})$ we obtain

\medskip

{\centerline{$F_2(p_K\restr\alpha_n,\alpha_{n+1}))\forces``\eta\in  T^*$.''}}

\medskip

Subclaim 2 is established.

To complete the induction establishing Claim 3, we use the Proper Iteration Lemma to take
$r_{n+1}\in P_{\alpha_{n+1}}$ such that $r_{n+1}\restr\alpha_n=r_n$ and
$r_{n+1}$ is $N$-generic and $r_n\forces``r_{n+1}\restr[\alpha_n,\alpha_{n+1})\leq
F_2(p_K\restr[\alpha_n,\alpha_{n+1}))$.''

Claim 3 is established.

Let $q'$ be a $P_\alpha$-name such that

\medskip

\centerline{$p\forces``q'=\bigcup\{r_n\restr[\alpha,\alpha_n)\,\colon\allowbreak n\in\omega\}$.''}

\medskip

In $V[G_{P_\alpha}]$, let $H^*=\bigcup T^*$ and let

\medskip

\centerline{$H=\{\nu\in H^*\,\colon\allowbreak(\forall n\in\omega)\allowbreak(\exists\eta\in H^*)\allowbreak
(\nu$ is comparable with $\eta)\}$.}

\medskip

 As in the proof of Lemma 7.12, we have that $H$ is a $z$-squeezed tree. By Claim 3 we have
that

\medskip

\centerline{$q'\forces``$for every $n\in\omega$ we have $\zeta_n\in[T^*]$ and $\zeta_n\restr n=\zeta\restr n$,}

\centerline{and therefore $\zeta\in[T^*]$,
and therefore $T\subseteq H$.''}

\medskip

The Theorem is established.

\proclaim Corollary 7.18.   Suppose $\langle P_\eta\,\colon\allowbreak\eta\leq\kappa\rangle$ is a countable support iteration
based on $\langle Q_\eta\,\colon\allowbreak\eta<\kappa\rangle$ and suppose that for every $\eta<\kappa$ we have that\/ {\rm
${\bf 1}\forces``Q_\eta$ is proper and
has the $P$-point property.''} Then $P_\kappa$ has
the $P$-point property.

Proof: By Theorem 7.17 with $\alpha=0$.

\section{On adding no Cohen reals}

In [12, Conclusion VI.2.13D(1)], Shelah states that a countable support iteration of proper forcings, each of which adds no Cohen
reals, either adds no Cohen reals or adds a dominating real.  However, according to Jakob Kellner,
Shelah has stated that this is an error, and the result holds only at limit stages.
In this section, we prove the limit case.

\proclaim Definition 8.1. A nowhere dense tree $T\subseteq{}^{<\omega}\omega$ is a non-empty tree such that
for every $\eta\in T$ there is some $\nu$ extending $\eta$ such that $\nu\notin T$. A perfect tree $T\subseteq{}^{<\omega}\omega$ is a 
non-empty tree such that
  for every $\eta\in T$, the set of successors of $\eta$ in $T$ is not linearly ordered.

\proclaim Lemma 8.2.  $P$ does not add any Cohen reals iff\/ {\rm ${\bf 1}\forces_P``(\forall f\in{}^\omega\omega)\allowbreak(\exists H\in V)\allowbreak(H$ is a nowhere dense perfect tree and $f\in[H])$.''}

Proof: This is a tautological consequence of the definition of Cohen real.

\proclaim Lemma 8.3.  Suppose ${\rm cf}(\kappa)=\omega$ and
$\langle P_\eta\,\colon\allowbreak\eta\leq\kappa\rangle$ is a countable
support forcing iteration. Suppose {\rm $(\forall\eta<\kappa)\allowbreak(P_\eta$  does not add any Cohen reals)} and\/
 {\rm
${\bf 1}\forces_{P_\kappa}``$for every countable $x\subseteq V$ there is a countable $y\in V$ such that
$x\subseteq y$.''}  Suppose $P_\kappa $ does not add any dominating reals. Then $P_\kappa$ does not add any Cohen reals.

Proof. Fix $\langle\alpha_n\,\colon\allowbreak n\in\omega\rangle$ cofinal in $\kappa$ with $\alpha_0=0$. Also in $V[G_{P_\kappa}]$ fix
$\zeta\in{}^\omega\omega$.
Use Lemma 3.1 to construct $\langle p_n\,\colon\allowbreak n\in\omega\rangle $ and $\langle\zeta_n\,\colon\allowbreak n\in\omega\rangle$
as there. In particular for every $n\in\omega$ we have $p_0\forces_{P_\kappa}``\zeta_n\restr n=\zeta\restr n$ and $\zeta_n\in
{}^\omega\omega\cap V[G_{P_{\alpha_n}}]$.''

Working in $V[G_{P_\kappa}]$ with $p_0\in G_{P_\kappa}$, use the fact that for every $n\in\omega$
we have that $P_{\alpha_n}$ does not add Cohen reals, let $\langle T_n\,\colon\allowbreak n\in\omega\rangle$ be a sequence of nowhere dense perfect trees such that
$(\forall n\in\omega)\allowbreak(T_n\in V$ and $\zeta_n\in[T_n])$.

 Let $B\in V$ be a countable set of nowhere dense perfect trees  such that for every $n\in\omega$ we have
$T_n\in B$.
Let $\langle S_n\,\colon\allowbreak n\in\omega\rangle\in V$ enumerate  $B$ with infinitely many repetitions such that $T_0=S_0$.

Build inductively $\langle S'_n\,\colon\allowbreak n\in\omega\rangle$ such that
$S'_{n+1}=S_{n+1}\cup S'_n$ and $S_0'=S_0$.

Define $h\in{}^\omega\omega$ by setting $h(k) $ equal to the least $m>k$ such that $T_{k}\subseteq S'_m$, for every $k\in \omega$.
Because $P_\kappa$ adds no dominating reals we may choose $g\in{}^\omega\omega\cap V$  and $A\subseteq \omega$ such that 
$A=\{n\in \omega\,\colon\allowbreak g(n)>h(n)\}$ and $A$ is infinite.

Choose $\langle k_i\,\colon\allowbreak i\in\omega\rangle\in V$ an increasing sequence of integers
as follows.
Let $k_0=0$. Given $k_n$, choose $k_{n+1}\geq{\rm max}(k_n+1,2)$ such that $(\forall\nu\in{}^{\leq k_n}k_n)\allowbreak
[(\exists\nu'\in{}^{k_{n+1}}\omega$ extending $\nu)\allowbreak(\forall i\leq k_n)\allowbreak(\nu'\notin S'_{g(i)})$ and
$(\forall i\leq k_n)\allowbreak(\exists \nu_1\in S'_{g(i)})\allowbreak(\exists \nu_2\in S'_{g(i)})\allowbreak(\nu_1$ and $\nu_2$ 
are distinct extensions of $\nu$ and
${\rm lh}(\nu_1)={\rm lh}(\nu_2)=k_{n+1})]$.

Let $T^0=\{\eta\in{}^{<\omega}\omega\,\colon\allowbreak(\exists s\in\omega)\allowbreak(\exists j\in\omega)\allowbreak
(k_{2s}\leq j<k_{2s+1}$ and $\eta\restr j\in S'_0$ and $\eta\in S'_{g(j)})\}$.

Let $T^1=\{\eta\in{}^{<\omega}\omega\,\colon\allowbreak(\exists s\in\omega)\allowbreak(\exists j\in\omega)\allowbreak
(k_{2s+1}\leq j<k_{2s+2}$ and $\eta\restr j\in S'_0$ and $\eta\in S'_{g(j)})\}$.

Claim 1: $T^0$ is a nowhere dense  tree.

Proof. Suppose $\eta\in T^0$. Choose $s$ and $j$ witnessing this. Also take $n\geq s$ so large that
$\eta\in{}^{\leq k_{2n}}k_{2n}$.

We  choose $\nu$ extending $\eta$ such that ${\rm lh}(\nu)=k_{2n+2}$ and
$(\forall i\leq k_{2n+1})\allowbreak(\nu\notin S'_{g(i)})$. In particular we have
$\nu\notin S'_0$.  We show that $\nu\notin T^0$.
So suppose, towards a contradiction, that $s'\in\omega$ and
$j'\in\omega$ and $k_{2s'}\leq j'<k_{2s'+1}$ and $\nu\restr j'\in S'_0$ and $\nu\in S'_{g(j')}$. Because $\nu\in S'_{g(j')}$ we know $j'\geq k_{2n+1}$.
Necessarily, then, $j'\geq k_{2n+2}$. But then $\nu=\nu\restr j'\in S'_0$.  This  contradiction establishes the 
Claim.

Claim 2. $T^0$ is a perfect tree.

Proof: Given $\eta\in T^0$, let $s\in\omega$ and $j\in\omega$ be witnesses. 

Case 1: ${\rm lh}(\eta)\geq j$.

Let $\nu$ and $\nu'$ be incomparable elements of $S'_{g(j)}$
extending $\eta$.  We have that $\nu$ and $\nu'$ are in $T^0$; this is witnessed by the integers $s$ and $j$.

Case 2: ${\rm lh}(\eta)<j$.

Take $\nu$ and $\nu'$ distinct extensions of $\eta$ such that $\nu\in S_0$ and $\nu'\in S_0$ and
${\rm lh}(\nu)={\rm lh}(\nu')=j$. We have $\nu\in S'_{g(j)}$ and $\nu'\in S'_{g(j)}$ because
$S_0\subseteq S'_{g(j)}$.  We have that $\nu$ and $\nu'$ are in $T^0$; this is witnessed by the integers $s$ and $j$.

Claim 3: $T^1$ is a nowhere dense perfect tree.

Proof: Similar to Claims 1 and 2.

Let $B_0=\bigcup\{[k_{2i},k_{2i+1})\,\colon\allowbreak i\in\omega\}$ and let $B_1=\bigcup\{[k_{2i+1},k_{2i+2})\,\colon\allowbreak i\in\omega\}$.

Claim 4:  $(\forall n\in A\cap B_0)\allowbreak(\zeta_{n}\in[T^0])$.
  
Proof: Given  $n\in A\cap B_0$ choose $s\in\omega$ such that
$k_{2s}\leq n< k_{2s+1}$. We have
$\zeta_{n}\in[T_{n}]\subseteq [S'_{h(n)}]\subseteq [S'_{g(n)}]$ and $\zeta_{n}\restr n\in S'_0$. Hence
$\zeta_{n}\in[T^0]$.  The Claim is established.

Claim 5:  $(\forall n\in A\cap B_1)\allowbreak(\zeta_{n}\in[T^1])$.

Proof: Similar to Claim 4.

We have that $T^0$ and $T^1$ are elements of $V$. Furthermore, if $A\cap B_0$ is infinite, we have by Claim 4  that 
for infinitely many $n$ we have $\zeta_{n}\in[T^0])$ and hence $\zeta\in[T^0]$.
Otherwise by Claim 5 it follows  that for infinitely many $n$ we have $\zeta_n\in[T^1]$ and hence $\zeta\in[T^1]$.  The Lemma is established.

\section{On not adding reals not belonging to any null sets of $V$}

In this section we give Shelah's proof that the property $``P$ does not add any real not belonging to any closed set of measure zero of the ground model'' is preserved at limit stages by countable support iterations of proper forcings assuming the iteration does not add dominating reals.

\proclaim Theorem 9.1.    Suppose $\langle P_\eta\,\colon\allowbreak\eta\leq\kappa\rangle$ is a countable support iteration based on $\langle Q_\eta\,\colon\allowbreak\eta<\kappa\rangle$ and suppose\/ $\kappa$ is a limit ordinal and\/ {\rm $(\forall\eta<\kappa)\allowbreak(P_\eta$ 
 does not add reals not in any closed measure zero set of $V$).}   Suppose also that $P_\kappa$ does not add any dominating reals.
Then $P_\kappa$ does not add any real not in any closed measure zero set of $V$.

Proof:  Repeat the proof of Theorem 8.3 with ``nowhere dense perfect tree'' replace by ``perfect tree with Lebesgue measure zero''
throughout.

\vfill\eject

\noindent{\bf Refernces}

\bigskip

[1] Todd Eisworth, CH and first countable, countably compact spaces, Topology Appl. 109, no. 1, 55-73 (2001)

\medskip

[2]  Goldstern, M., Tools for Your Forcing Construction, Set theory of the reals (Haim Judah, editor), Israel Mathematical Conference Proceedings, vol. 6, American Mathematical Society, pp. 305--360. (1993)

\medskip

[3]  Goldstern, M. and J. Kellner, New reals: Can live with them, can live without them, Math. Log. Quart. 52, No. 2, pp.~115--124 (2006)

 \medskip

[4] Kellner J., and S. Shelah, Preserving preservation, Journal of Symbolic Logic, vol 70, pp.~914--945 (2005)

\medskip

[5] Schlindwein, C., Consistency of Suslin's hypothesis, a non-special Aronszajn tree, and GCH, Journal of Symbolic Logic, vol.~59, pp.~1--29, 1994

\medskip

[6] Schlindwein, C., Suslin's hypothesis does not imply stationary antichains, Annals of Pure and Applied Logic, vol.~64, pp.~153--167, 1993

\medskip

[7] Schlindwein, C., Special non-special $\aleph_1$-trees, Set theory and its Applications,  J.~Steprans and S.~Watson (eds.), Lecture
Notes
in Mathematics, vol.~1401, Springer-Verlag, 1989

\medskip

[8] Schlindwein, C., A short proof of the preservation of the $\omega^\omega$-bounding property, Mathematical Logic Quarterly, vol.~50, pp.~19--32, 2004

\medskip

[9] Schlindwein, C., Shelah's work on non-semi-proper iterations, II, Journal of Symbolic Logic, 2001

\medskip

[10] Schlindwein, C., SH plus CH does not imply stationary antichains, Annals of Pure and Applied Logic, vol.~124, pp.~233--265, 2003

\medskip

[11] Schlindwein, C., Shelah's work on non-semi-proper iterations, I, Archive for Mathematical Logic.

\medskip

[12] Shelah, S., {\bf Proper and Improper Forcing}, Perspectives in Mathematical Logic, Springer, Berlin, 1998

\end{document}